\documentclass[10pt]{article}
\usepackage{amssymb}
\usepackage{amsmath}
\usepackage{amsthm}
\theoremstyle{plain}
\newtheorem{theorem}{Theorem}[section]
\newtheorem{proposition}[theorem]{Proposition}
\newtheorem{lemma}[theorem]{Lemma}
\newtheorem{corollary}[theorem]{Corollary}
\theoremstyle{definition}

\theoremstyle{remark}
\newtheorem{remark}[theorem]{Remark}

\newcommand{\N}{\mathbb{N}}

\newcommand{\Z}{\mathbb{Z}}
\newcommand{\R}{\mathbb{R}}
\newcommand{\C}{\mathbb{C}}

\newcommand{\T}{\mathbb{T}}
\newcommand{\oalpha}{\overline{\alpha}}

\newcommand{\Zj}{{\cal Z}}

\newcommand{\Oi}{{\cal O}_\infty}
\newcommand{\heta}{\hat{\eta}}
\newcommand{\hpsi}{\hat{\psi}}
\newcommand{\hxi}{\hat{\xi}}
\newcommand{\td}{\widetilde{d}}
\newcommand{\tp}{\widetilde{p}}

\newcommand{\tI}{\widetilde{I}}
\newcommand{\tiota}{\widetilde{\iota}}
\newcommand{\tX}{\widetilde{X}}

\newcommand{\tphi}{\widetilde{\varphi}}

\DeclareMathOperator{\id}{id}

\DeclareMathOperator{\Aut}{Aut}
\DeclareMathOperator{\Ad}{Ad}
\DeclareMathOperator{\WInn}{WInn}

\DeclareMathOperator{\rank}{rank}

\DeclareMathOperator{\diag}{diag}

\DeclareMathOperator{\Tr}{Tr}

\newcommand{\bprf}{\noindent{\it Proof.\ }}

\newcommand{\bprfm}{\noindent{\it Proof of Theorem 1.2. \quad}}
\newcommand{\eprf}{\hspace*{\fill} \rule{1.6mm}{3.2mm} \vspace{1.6mm}}
\newcommand{\benu}{\begin{enumerate}\renewcommand{\labelenumi}{{\rm (\roman{enumi})}}\renewcommand{\itemsep}{0pt}}
\newcommand{\eenu}{\end{enumerate}}
\begin{document}

\title{A generalization of the Jiang-Su construction }
\author{Yasuhiko Sato \\}
\date{\small  Department of Mathematics, Hokkaido University, \\ Sapporo 060-0810, Japan \\ e-mail : s053018@math.sci.hokudai.ac.jp}

\maketitle

\begin{abstract}
 Let $G$ be a countable abelian group. We construct a unital simple projectionless $C^*$-algebra $A$ with a unique tracial state, that satisfies $(K_0(A), [1_A]) \cong (\Z, 1) $, $K_1(A) \cong G$,  $A\otimes \Zj \cong A$, and that is obtained as the inductive limit $C^*$-algebra of a sequence of dimension drop algebras of a specific form. This construction is based on the construction of the Jiang-Su algebra. By this construction, we show a certain conjugacy result for aperiodic automorphisms of these projectionless $C^*$-algebras. We also show that an automorphism of this projectoinless $C^*$-algebra has a certain aperiodicity up to the weakly inner automorphisms in the tracial representation if and only if it has a kind of Rohlin property, which leads to the Rohlin property after taking tensor product of certain $C^*$-algebras of real rank zero.   



\end{abstract}

\section{Introduction}\label{INTRO}
 In \cite{JS}, Jiang and Su have concretely constructed the projectionless $C^*$-algebra $\Zj$ whose $K_0$-group is $(\Z, \Z_+, 1)$ and $K_1$-group is $\{0\}$ which is now called the Jiang-Su algebra and they have shown that any automorphism of the Jiang-Su algebra is approximately inner. For the classification program of finite $C^*$-algebras started by G. A. Elliott , the Jiang-Su algebra behaves as $\Oi$ in the classification theory for Kirchberg algebras. Then the Jiang-Su algebra was studied by many people (e.g., \cite{DW}, \cite{DPT}, \cite{RW}, \cite{TW}, and \cite{W} ).     

In this paper elaborating Jiang-Su's method we give a unified way of constructing unital projectionless $C^*$-algebras with arbitrary $K_1$-group.  And, as an application of this construction, we can prove the following corollary as remarked in Remark \ref{remarkcor}.
\begin{corollary}\label{maincorollary}
Let $A$ be the unital simple projectionless $C^*$-algebra that is obtained in  Theorem \ref{maintheorem} or \ref{maintheorem2}, and $\tau$ the unique tracial state of $A$. Suppose that $\alpha $ and $\beta$ are automorphisms of $A$ which are aperiodic in the quotient $\Aut(A)/\WInn(A)$, then for any $\varepsilon >0$ there exists an approximately inner automorphism $\sigma \in \Aut (A)$ and $\gamma \in \WInn (A) $ with $W \in U(\pi_{\tau}(A)''
)$ such that $\pi_{\tau} \circ \gamma =\Ad W \circ \pi_{\tau}$, $\|W-1 \|_2< \varepsilon $, and
\[\alpha= \gamma \circ \sigma \circ \beta \circ \sigma^{-1}. \]
\end{corollary}
Here $\WInn (A)$ is the subgroup of the automorphism group of $A$ defined by
\[ \WInn (A) = \{ \alpha \in \Aut (A) ; \pi_{\tau} \circ \alpha = \Ad V\circ \pi_{\tau},\ V \in U(\pi_{\tau} (A)'')\},\]
 and $\pi_{\tau} $ is the GNS-representation associated with $\tau$. 
This corollary is a straightforward generalization of Theorem 2.6 in \cite{Sat}. The aperiodicity in the assumption of this corollary will be characterized by the following theorem.
\begin{theorem}\label{maintheorem3}
Let $A$ and $\tau$ be as in the above corollary and  $\alpha $ is an automorphism of $A$. 
Then $[\alpha] \in \Aut (A) / \WInn(A) $
 is aperiodic
if and only if $\alpha$ satisfies the following condition that:
for any $k \in \N$ there exists a central sequence $(f_n)_n $ of positive elements in $A$ (i.e.,  $(f_n)_n \in (A_{\infty})_+^1$) such that 
\[ (f_n)_n \cdot \alpha^j ((f_n)_n) = 0,\ j=1,2,...,k-1,\quad  \tau (1_A - \sum_{j=0}^{k-1} \alpha^{j} (f_n) ) \rightarrow 0 ,\] 
\end{theorem}

This theorem is motivated by a similar result for UHF-algebras by Kishimoto \cite{Kis}. The condition expressed in terms of positive elements in this theorem, has given the result in \cite{Sat} that: when $A$ is the Jiang-Su algebra, $B$ is a unital simple TAF-algebra satisfying a certain condition in \cite{Sat} , $\beta$ is  an automorphism of $B$, and $\alpha\in \Aut (A)$ has the property in the above theorem, this condition implies that the automorphism $\alpha\otimes \beta$ of $A\otimes B$ is asymptotically unitarily equivalent to an automorphism with the Rohlin property. 

The building blocks of the unital simple projectionless $C^*$-algebras constructed in Theorem \ref{maintheorem} or \ref{maintheorem2} will be  
dimension drop algebras on 
one-dimensional compact Hausdorff graphs with `dimension drop' condition at each end point. Give a finitely generated abelian group $G$, the base space is uniquely determined as a one-point union of copies of torus $\T$ and copies of the interval $[0,2\pi]$ with $0 \in [0, 2\pi]$ the common point, where the number of copies of $\T$ corresponds to the rank of the free part of $G$ and the number of copies of $[0,2\pi]$ corresponds to the number of directly indecomposable summands of the torsion subgroup of $G$. Regulating matrix algebras on the end points $2\pi \in [0, 2\pi]$, we obtain a unital projectionless dimension drop algebra that has the finite cyclic groups as $K_1$-group.
In \cite{Ell},  Elliott have already obtained the projectionless $C^*$-algebra whose  $K_0$-group is $(\Z,\Z_+,1)$ and these $K_1$-group is a given countable abelian group, and many people have studied subhomogeneous $C^*$-algebras on one dimensional compact Hausdorff spaces of similar type (see e.g. \cite{Li}, \cite{My}, \cite{Thomsen}). 
However the author could not find a type of unital projectionless dimension drop algebras suitable to prove Corollary \ref{maincorollary}. Recently, H. Lin showed the classification result of the inductive limit $C^*$-algebras of dimension drop algebras of general type, in \cite{Lin}.

The contents of this paper are as follows. In Section \ref{Sec2}, for a given finitely generated abelian group $G$ we construct a projectionless dimension drop algebra whose $K_1$-group is $G$. In Section \ref{Sec3}, for a given countable abelian group $G$ we construct a unital simple projectionless $C^*$-algebra whose $K_1$-group is $G$ that is the inductive limit $C^*$-algebra of a sequence of dimension drop algebras. In Section \ref{Sec4}, we show $\Zj$-stability of these projectionless $C^*$-algebras, which is a key property in the classification theory according to  \cite{ET} and \cite{W}. In Section \ref{Sec5}, we prove Corollary \ref{maincorollary} and Theorem \ref{maintheorem3}. In order to show this corollary we need to show a UHF-embeddability, and in order to show this UHF-embeddability we also need to consider another one-dimensional homogeneous $C^*$-algebra which  naturally contains the dimension drop algebra mentioned above. 

Concluding this section, we prepare some notations. When $A$ is a $C^*$-algebra, we denote by $A^1$ the set $\{a\in A : \| a \| \leq 1 \}$, $A_{\rm sa}$ the set of self-adjoint elements of $A$, $U(A)$ the unitary group of $A$, $T(A)$ the  tracial state space of $A$, $A^{\infty}$ the quotient $\ell^{\infty}(\N,A)/c_0(A)$, and $A_{\infty}$
 the relative commutant $A^{\infty} \cap A'.$
 We use $M_n$ to denote the $C^*$-algebra of $n\times n$ matrices with complex entries and we denote by $\Tr_n $  the normalized trace on $M_n$.
In this paper we denote by $p_i$ the $i$th-prime number and $(m,n)$ the greatest common divisor of $m$ and $n$. On a metric space $(X,d)$, we say that $Y \subset X$ is {\it $\varepsilon$-dense }in $X$ if for $y \in Y$ there is $x \in X $ such that $d(x,y) \leq \varepsilon $.

\section{Dimension drop algebras}\label{Sec2}
First we shall prepare some symbols and define a dimension drop algebra for a finitely generated abelian group. 
For a given countable abelian group $G$ we will express $G$ as the inductive limit of an increasing sequence of finitely generated abelian group $G_n$, $n \in \N$, and construct a unital simple projectionless $C^*$-algebra whose $K_1$-group is $G$ as the inductive limit $C^*$-algebra of dimension drop algebras for $G_n$, $n\in \N$. Then we denote by $G_n$ a finitely generated abelian group, by $G_{n,0}$ the free part of $G_n$, and by $G_{n,1}$ the torsion subgroup of $G_n$. Set $r_0^{(n)} = \rank (G_{n,0}) \in \Z_{+}$. For $i \in \N$ with $p_i |\  | G_{n,1}| ,$ let $G_{n,i}$ be the subgroup of $G_n$ generated by $\{ g\in G_n ; o(g)= p_i^{n}, n \in \N  \}$, where $o(g)$ means the order of $g$, and let $r_i^{(n)} \in \N$ and $d_{i,j}^{(n)} \in \N,$ $j=1,2,...,r_i^{(n)}$ be such that 
\[ G_{n,i} \cong \bigoplus_{j=1}^{r_i^{(n)}} \Z_{d_{i,j}^{(n)}}, \quad d_{i,j}^{(n)} = p_i ^{k_{i,j}^{(n)}},\] 
where $\Z_{d_{i,j}^{(n)}}=\Z/ d_{i,j}^{(n)} \Z$, $k_{i,j}^{(n)} \in \N$ and $k_{i,j}^{(n)} \leq k_{i,j+1}^{(n)}$.
Therefore we identify  
\[ G_n \cong G_{n,0} \oplus G_{n,1} \cong \Z ^{r_0^{(n)}} \oplus \bigoplus_{p_i |\ |G_{n,1}|} \ \bigoplus_{j=1}^{r_i^{(n)}} \Z_{d_{i,j}^{(n)}} . \]
We define canonical generators $g_{i,j}^{(n)}$ of $G_n$ by
\[g_{0,j}^{(n)}= (0\oplus \cdots \oplus 0 \oplus 1_{j} \oplus 0 \oplus \cdots \oplus 0) \oplus 0_{G_{n,1}},\  j =1,2,..., r_0^{(n)}, \]
\[g_{i,j}^{(n)} = 0_{G_{n,0}} \oplus (0 \oplus \cdots \oplus 0 \oplus 1_{i,j} \oplus 0 \oplus \cdots \oplus 0),\ j = 1,2,..., r_i^{(n)}.\]

For a finitely generated abelian group $G_n$, set 
\[S^{(n)} = \{ (0,0), (0,1)\} \cup  \{ (i,j) ;\ i\in \N \ {\rm with}\ p_i|\ |G_{n,1}|,\ j=1,2,...,r_i^{(n)} \}. \] 
 Let $I_{i,j}^{(n)}$ be a copy of the interval $[0, 2\pi]$ for each $(i,j) \in S^{(n)}$ and denote by $\iota_{i,j}^{(n)}$ the identity map from $[0,2\pi]$ onto $I_{i,j}^{(n)}$, and let $T_{j}^{(n)}$ be a copy of the torus $\T = \{ e^{it}; t\in [0, 2 \pi] \} $ for each $j=1,2,...,r_0^{(n)} $ and denote by $\tau_j^{(n)}$ the identity map from $\T$ onto $T_j^{(n)}$. Set $z(t)=e^{it} \in \T $, $t\in \R$. We define a one-point union $X_n$ by 
\[ X_n = I_{0,0}^{(n)} \vee I_{0,1}^{(n)} \vee \bigvee_{j=1}^{r_0^{(n)}} T_j^{(n)} \vee \bigvee_{p_i |\ |G_{n,1}|}\ \bigvee_{j=1}^{r_i^{(n)}} I_{i,j}^{(n)} , \]
 identifying the base point of $X_n$ with $\iota_{i,j}^{(n)}(0)$, $(i,j)\in S^{(n)} $ and $\tau_j^{(n)} (z(0))$, $j=1,2,...,r_0^{(n)}$.
Set $c_n= \tau_j^{(n)} (z(0)) = \iota_{i,j}^{(n)} (0) \in X_n $ and $a_{i,j}^{(n)} = \iota_{i,j}^{(n)} (2 \pi) \in X_n$ for $(i,j) \in S^{(n)}$.

For a finitely generated abelian group $G_n$ and relatively prime natural numbers $\tp_0$ and $\tp_1 $ such that  $ d_{i,j}^{(n)} \mid \tp_0 \cdot \tp_1 $ for any $(i,j)\in S^{(n)}$ with $i\geq 1$, we set $\natural (i) = 1-i$ for $i=0,1$ and $\natural (\mu, \nu) = \natural (i)$ for $(\mu,\nu) \in S^{(n)}$ with $\mu \geq 1$ and $p_{\mu}|\tp_{i}$, and we define a $C^*$-algebra $A(G_n, \tp_0, \tp_1)$ by
\begin{eqnarray}
A(G_n, \tp_0, \tp_1) &=& \{ f \in C(X_n)\otimes M_{\td}\ ; f(a_{0,i}^{(n)}) \in \epsilon_{0,i}^{(n)} (M_{\tp_0}),\ i=0,1, \nonumber \\
& & f(a_{i,j}^{(n)})\in \epsilon_{i,j}^{(n)} (M_{\td/d_{i,j}^{(n)}}),\ (i,j)\in S^{(n)} \ {\rm with}\ i \geq 1 \},\nonumber 
\end{eqnarray}
where $ \td = \tp_0 \cdot \tp_1 $ ( and $d_{i,j}^{(n)} = p_i^{k_{i,j}^{(n)}}$) and 
$\epsilon_{0,i}^{(n)} $ and $\epsilon_{i,j}^{(n)}$ are unital embeddings from $M_{\tp_i}$ and $M_{\td /d_{i,j}^{(n)}}$ onto $M_{\td}$. It is not so hard to show that the $C^*$-algebra $A(G_n, \tp_0, \tp_1)$ is independent of choices of $\epsilon_{i,j}^{(n)}$, $(i,j)\in S^{(n)}$.
In this paper we say that $A(G_n, \tp_0, \tp_1) $ is a {\it dimension drop algebra for $G_n$, $\tp_0$, and $\tp_1$}.
 In particular we denote by $I(\tp_0, \tp_1)$ the dimension drop algebra for $\{0\},$ $\tp_0,$ and $\tp_1$  with ($\tp_0, \tp_1$) =1, (i.e., 
\[ I( \tp_0, \tp_1) = \{ f\in C( I_{0,0}^{(n)} \cup I_{0,1}^{(n)}) \otimes M_{\td} \ ; f (a_{0,0}^{(n)}) \in M_{\tp_0} \otimes 1_{\tp_1} , \ f(a_{0,1}^{(n)}) \in 1_{\tp_0} \otimes M_{\tp_1} \}), \]
that was defined in \cite{JS}.

Let $e_{0,0}$ be a minimal projection of $M_{\td}$. We define $f_{0,j}^{(n)} \in A(G_n, \tp_0, \tp_1) $ for  $j=1,2,...,r_0^{(n)}$ and  $f_{i,j}^{(n)} \in A(G_n, \tp_0, \tp_1)$ for $(i,j)\in S^{(n)}$ with $i\geq 1$ by 
\begin{eqnarray}
f_{0,j}^{(n)}(x)=\left\{ \begin{array}{ll}
\exp (ite_{0,0}) ,\quad & x=\tau_j^{(n)}(z(t))\in T_j^{(n)},\ t \in [0,2\pi], \\
1_{\td}, \quad  & {\rm otherwise}, \\
\end{array} \right.\nonumber
\end{eqnarray}   

\begin{eqnarray}
f_{i,j}^{(n)}(x)=\left\{ \begin{array}{ll}
\exp (ite_{0,0}) ,\quad & x= \iota_{i,j}^{(n)}(t)\in I_{i,j}^{(n)},\ t \in [0,2\pi], \\
1_{\td}, \quad  & {\rm otherwise}. \\
\end{array} \right.\nonumber
\end{eqnarray} 

\begin{proposition}\label{Prop1}
Let $G_n$, $\tp_i$, $i=0,1$, and $A(G_n,\tp_0, \tp_1)$ be as above. 
Then it follows that $(K_0(A(G_n,\tp_0,\tp_1)), [1]) \cong (\Z ,1)$ and that 
\[ K_1 (A(G_n,\tp_0,\tp_1)) \cong G_n ,\]
under the identification determined by  
\[ [f_{0,j}^{(n)}]_{K_1(A(G_n,\tp_0,\tp_1)} \cong g_{0,j}^{(n)},\ j=1,2,..., r_0^{(n)}, \quad [f_{i,j}^{(n)} ]_{K_1 (A(G_n,\tp_0,\tp_1)} \cong g_{i,j}^{(n)}, \]
for $(i,j) \in S^{(n)} $ with $i \geq 1$.
\end{proposition}  
\bprf
Set $A_n= A(G_n, \tp_0, \tp_1)$ and $I = I ( \tp_0, \tp_1)$, then
we have that $K_1 ( I ) \cong  0  $ from $(\tp_0, \tp_1) =1 $.
Set  
\begin{eqnarray}
B_{n,0}&=& \bigoplus_{(i,j)\in S^{(n)},\ i\geq 1} M_{\td/d_{i,j}^{(n)}}, \nonumber \\
X_{n,0} &=& \bigcup_{j=1}^{r_0^{(n)}} \tau_j^{(n)} (z( (0, 2 \pi))) \cup \bigcup_{(i,j)\in S^{(n)},\ i\geq 1} \iota_{i,j}^{(n)} ((0,2 \pi)) \subset X_n, \nonumber \\C_{n,0} &=& C_0 (X_{n,0}) \otimes M_{\td} . \nonumber
\end{eqnarray}
Then we have the exact sequence 
\[ 0 \longrightarrow C_{n,0} \longrightarrow\hspace{-1.25em}{^{\rho}}\ \  A_n \longrightarrow\hspace{-1.25em}{^{\sigma}}\ \ I \oplus B_{n,0} \longrightarrow 0, \] 
which induces the exact sequence on $K$-groups 
\[K_0 (I \oplus B_{n,0}) \longrightarrow\hspace{-1.25em}{^{\delta_0}}\ \  K_1((C_{n,0})^{\sim})\longrightarrow K_1(A_n) \longrightarrow 0,  \]
where $\delta_0$ is the exponential map.

By $ (\tp_0,\tp_1) =1 $, it follows that $(K_0(I),[1_{I}]_{K_0(I)}) \cong (\Z, 1) $. Let $e_{i,j}$ be a minimal projection of $M_{\td/d_{i,j}^{(n)}} \subset B_{n,0} $. It follows that 
$ K_0(B_{n,0}) \cong \bigoplus_{(i,j)\in S^{(n)}, i\geq 1}$ $ \Z $ under the identification defined by 
\[ [e_{i,j}]_{K_0 (B_{n,0})} \cong 0\oplus\cdots\oplus 0 \oplus 1_{i,j} \oplus0\oplus \cdots \oplus  0 ,\] 
and it follows that $ \displaystyle K_1((C_{n,0})^{\sim}) \cong \bigoplus_{j=1}^{r_0^{(n)}}  \Z \oplus  \bigoplus_{(i,j)\in S^{(n)},\ i\geq1} \Z $ by the identification  defined by
\[ [f_{0,j}^{(n)}]_{K_1 ((C_{n,0})^{\sim})} \cong (0 \oplus \cdots \oplus 0 \oplus 1_j \oplus 0 \oplus \cdots \oplus 0 ) \oplus 0_{\bigoplus_{i,j} \Z},\]
\[ [f_{i,j}^{(n)} ]_{K_1( (C_{n,0})^{\sim})} \cong 0_{\bigoplus_j \Z} \oplus (0\oplus 
\cdots \oplus 0 \oplus 1_{i,j} \oplus 0 \oplus \cdots \oplus 0 ) .\]
 
Define $h \in (A_n)_{\rm sa}$ by 
\begin{eqnarray}
h(x)=\left\{ \begin{array}{ll}
 (1-\frac{t}{2\pi})1_{\td} ,\quad & x = \iota_{i,j}^{(n)} (t) \in I_{i,j}^{(n)},\ i\geq 1, \\
1_{\td}, \quad  & {\rm otherwise}. \\
\end{array} \right.\nonumber
\end{eqnarray}
Because of the definition of $\delta_0$ and $\sigma (h) = 1_{I}\oplus 0_{B_{n,0}}$, we have that
\[ \delta_0 ([1_{I}\oplus 0_{B_{n,0}} ]_{K_0 (I \oplus B_{n,0})}) =- [\exp (2\pi i h)]_{K_1 ((C_{n,0})^{\sim})} \cong 0_{\bigoplus_j \Z} \oplus (\td \oplus \cdots \oplus \td) ,\]
under the above identification. Define $h_{i,j}^{(n)} \in (A_n)_{\rm sa}$ for $(i,j) \in S^{(n)}$ with $i \geq 1$ by
\begin{eqnarray}
h_{i,j}^{(n)} (x)=\left\{ \begin{array}{ll}
\frac{t}{2\pi} \epsilon_{i,j}^{(n)} (e_{i,j}) ,\quad & x = \iota_{i,j}^{(n)} (t) \in I_{i,j}^{(n)},\ t \in[0,2\pi], \\
0_{\td}, \quad  & {\rm otherwise}. \\
\end{array} \right.\nonumber
\end{eqnarray} 
Since $\sigma (h_{i,j}^{(n)} ) =0_{I} \oplus e_{i,j}$, it follows that 
\begin{eqnarray}
\delta_0([0_{I} \oplus e_{i,j}]) &=& -[ \exp (2\pi i h_{i,j}^{(n)}) ]_{K_1((C_{n,0})^{\sim})} = -d_{i,j}^{(n)} [ f_{i,j}^{(n)} ] \nonumber \\
&=& 0_{\bigoplus_j \Z} \oplus (0\oplus \cdots \oplus 0 \oplus - d_{i,j}^{(n)} \oplus 0 \oplus \cdots \oplus 0 ) . \nonumber
\end{eqnarray}
Therefore it follows that 
\begin{eqnarray}
K_1(A_n) &\cong & K_1( (C_{n,0})^{\sim})/ \delta_0 ( K_0 (I_{n} \oplus B_{n,0})) \nonumber \\
&\cong & \bigoplus_{j=1}^{r_0^{(n)}} \Z \oplus \bigoplus_{(i,j)\in S^{(n)},\ i\geq 1} \Z_{d_{i,j}^{(n)}}, \nonumber 
\end{eqnarray}
and that $[ f_{0,j}^{(n)}]_{K_1(A_n)}$ and $[f_{i,j}^{(n)} ]_{K_1(A_n)}$ correspond to $g_{0,j}^{(n)}$ and $g_{i,j}^{(n)}$ because of the above identification $K_1((C_{n,0})^{\sim}) \cong \bigoplus_{j} \Z \oplus  \bigoplus_{i,j} \Z. $

Since $K_0 (C_{n,0}) \cong 0 $, we have the exact sequence 
\[0 \longrightarrow K_0(A_n) \longrightarrow\hspace{-1.50em}{^{\sigma_*}}\ \  K_0 (I \oplus B_{n,0}) \longrightarrow\hspace{-1.25em}{^{\delta_0}}\ \  K_1((C_{n,0})^{\sim}).  \]
Then we have that $(K_0(A_n), [1_{A_n}]) \cong (\ker (\delta_0), [ \sigma (1_{A_n})])$. Because $[\sigma(1_{A_n})] = 1 \oplus (\td / d_{i,j}^{(n)})_{i,j}$ is the generator of $\ker(\delta_0)$ we conclude that $(\ker(\delta_0), [ \sigma (1_{A_n})]) \cong (\Z, 1).$ 
\eprf

The following definitions are necessary to show the UHF-embeddability in Proposition \ref{prop:UHF} of the projectionless $C^*$-algebra which will be obtained in Theorem \ref{maintheorem}. For a finitely generated abelian group $G_n$, we denote by $\tI_{0,j}^{(n)}$, $j=1,2,...,r_0^{(n)}$, $ r_0^{(n)}$ copies of the interval $[0, \pi]$, by $\tI_{1,j}^{(n)}$, $r_0^{(n)}$ copies of the interval $[-\pi, 0]$, by $\tiota_{0,j}^{(n)}$ the identity map $[0, \pi] \rightarrow \tI_{0,j}^{(n)}$, and by $\tiota_{1,j}^{(n)}$ the identity map $[-\pi, 0] \rightarrow \tI_{1,j}^{(n)}$. We define a one-point union $\tX_n$ by 
\[ \tX_n = I_{0,0}^{(n)} \vee I_{0,1}^{(n)} \vee \bigvee_{j=1}^{r_0^{(n)}} (\tI_{0,j}^{(n)} \vee \tI_{1,j}^{(n)})  \vee \bigvee_{p_i |\ |G_{n,1}|}\ \bigvee_{j=1}^{r_i^{(n)}} I_{i,j}^{(n)} , \]
identifying the base point of $\tX_n$ with $\iota_{i,j}^{(n)}(0)$ and $\tiota_{i,j}^{(n)}(0)$. We set $b_{i,j}^{(n)}= \tiota_{i,j}^{(n)}((-1)^i\pi)$ for $i=0,1$ and $j=1,2,...,r_0^{(n)}$, and we identify $A(G_n, \tp_0, \tp_1)$ with the $C^*$-subalgebra $\{f\in C(\tX_n)\otimes M_{\td};$ $f(b_{0,j}^{(n)})=\epsilon_j^{(n)} (f(b_{1,j}^{(n)})),$ $f(a_{0,i}^{(n)})\in \epsilon_{0,i}^{(n)} (M_{\tp_i}),$ $f(a_{i,j}^{(n)})\in \epsilon_{i,j}^{(n)} (M_{\td/d_{i,j}^{(n)}})\}$,
where $\epsilon_j^{(n)}$, $j=1,2,...,r_0^{(n)}$ are fixed automorphisms of $M_{\td}$, and we denote the set of end points of $\tX_n$ by  
\[E_n=\{ a_{i,j}^{(n)} ; (i,j) \in S^{(n)}\}\cup \{b_{i,j}^{(n)} ;\ i=0,1,\  j=1,2,...,r_0^{(n)}\}.\] 

Given $G_n$ and $A(G_n, \tp_0, \tp_1)$ as in Proposition \ref{Prop1} let $\gamma_n$ be an injective homomorphism from $G_n$ into another finitely generated abelian group $G_{n+1}$. We shall  show how to lift $\gamma_n$ to an embedding of $A(G_n, \tp_0, \tp_1)$ into another dimension drop algebra $A(G_{n+1}, \tp_0', \tp_1')$.
\begin{proposition}\label{Prop2}
Let $\gamma_n$ be an injective group homomorphism from $G_n$ into $G_{n+1}$.
 Then there exist relatively prime numbers $\tp_0'$ and $\tp_1'$ such that $d_{i',j'}^{(n+1)} |\ \tp_0' \tp_1'=\td'$ for any $(i',j')\in S^{(n+1)}$ with $i'\geq 1$, automorphisms $\{ \epsilon_{j'}^{(n+1)}\in \Aut (M_{\td'});$ $j=1,2,...,r_0^{(n+1)}\}$, unital embeddings $\{\epsilon_{i',j'}^{(n+1)}:M_{\td'/d_{i',j'}^{(n+1)}}\hookrightarrow M_{\td'};$ $(i',j')\in S^{(n+1)}\}$, where $d_{0,i}^{(n+1)}=\tp'_{\natural (i)}$ for $i=0,1$, and a unital embedding $\psi_n : C(\tX_n)\otimes M_{\td} \hookrightarrow C(\tX_{n+1})\otimes M_{\td'}$ satisfying the following conditions: 
\begin{description}\label{conditions:psi}
\item[(0,0)]\label{(0,0)} $\psi_n(A(G_n,\tp_0,\tp_1)) \subset A(G_{n+1}, \tp_0', \tp_1')$ and $(\psi_n|_{A(G_n,\tp_0,\tp_1)})_* \cong \gamma_n $ under the identification that was defined in Proposition \ref{Prop1}, where $A(G_{n+1},\tp_0',\tp_1')$ are defined by $\epsilon_{j'}^{(n+1)}$ and $\epsilon_{i',j'}^{(n+1)}$,  and $(\psi_n|_{A_n})_*$ is the induced map on $K_1$-groups, 
\item[(0,1)]\label{(0,1)} for $\varepsilon > 0$ and $f\in C(\tX_n) \otimes M_{\td}$, if $\psi_n (f) ^{-1} (\{0 \})$ is $ \varepsilon$-dense in $\tX_{n+1}$ then  $f^{-1} (\{0 \})$ is $ \varepsilon$-dense in $\tX_{n} $,
\item[(0,2)]\label{(0,2)} for any $\varepsilon >0 $ and $f\in C(\tX_n) \otimes M_{\td}$, 
\[ \max \{ |\Tr_{\td'}(\psi_{n}(f) (x)- \psi_{n} (f)(y))| ;\ x,y  \in \tX_{n+1},\ d(x,y) \leq \varepsilon \}\] 
\[ \leq \max \{ |\Tr_{\td} (f(x)-f(y))| ;\ x,y \in \tX_n,\ d(x,y) \leq  \varepsilon\},\]
where $d$ is the canonical metric on $\tX_n$,
\item[(0,3)]\label{(0,3)} $\psi_n (\chi (E_n) \otimes 1_{\td}) \leq \chi (E_{n+1}) \otimes 1_{\td'}$. 
\end{description} 

\end{proposition}\bprf
Suppose that $G_n$ and $G_{n+1}$ are defined by $(r_i^{(n)}, d_{i,j}^{(n)})$ and $(r_i^{(n+1)}, d_{i,j}^{(n+1)})$ as in the beginning of this section and let $\gamma_n: G_n \rightarrow G_{n+1}$ be an injective group homomorphism.
For the free part of $G_n$ and $j=1,2,..., r_0^{(n)}$, there exist $\mu_{j',j}\in \Z$, $j'=1,2,...,r_0^{(n+1)}$ and $\nu_{i',j',0,j} \in \Z_+$, $(i',j')\in S^{(n+1)}$ with $i' \geq 1 $ such that 
\[ \gamma_n(g_{0,j}^{(n)}) = \sum_{j'=1}^{r_0^{(n+1)}} \mu_{j',j} g_{0,j'}^{(n+1)} + \sum_{(i',j')\in S^{(n+1)},\ i'\geq 1} \nu_{i',j',0,j} g_{i',j'}^{(n+1)},\]
and that $0 \leq \nu_{i',j',0,j} < d_{i',j'}^{(n+1)}(=p_{i'}^{k_{i',j'}^{(n+1)}})$. 
 For the torsion group $G_{n,1}$ and  $(i,j) \in S^{(n)}$ with $i\geq1$, because of $\gamma_n (G_{n,p_i}) \subset G_{n+1,p_i}$, there exist $\nu_{j',i,j}\in \Z$, $j'=1,2,...,r_{i'}^{(n+1)}$ such that 
\[ \gamma_n(g_{i,j}^{(n)}) = \sum_{j'=1}^{r_i^{(n+1)}} \nu_{j',i,j} g_{i,j'}^{(n+1)},\] and that $0 \leq \nu_{j',i,j} < d_{i,j'}^{(n+1)}$. 
When $ p_{i'} \nmid \ |G_{n,1}|  $ we define $r_{i'}^{(n)} = 0,$ $d_{i',j}^{(n)} = 1$ and $\nu_{j',i',j}=0$. Remark that 
$ d_{i',j'}^{(n+1)} |  d_{i',j}^{(n)} \nu_{j',i',j} $ (by $0 = \gamma_n(d_{i',j}^{(n)}g_{i',j}^{(n)}) = \sum d_{i',j}^{(n)} \nu_{j',i',j} g_{i',j'}^{(n+1)}$).  

Let $P_0$ and $P_1$ be large natural numbers  such that
\[ P_0P_1 > \max (\{ \sum_{j=1}^{r_0^{(n)}} |\mu_{j'j}| ;\ j'=1,2,...,r_0^{(n+1)}\} \cup \{(r_0^{(n)}+r_{i'}^{(n)})d_{i',j'}^{(n+1)} ;\ i', j'\} ) ,\]
\[ (P_0 \tp_0, P_1 \tp_1) =1 , \quad {\rm and} \ d_{i',j'}^{(n+1)}\ | P_0 P_1,\]
for any $(i', j')\in S^{(n+1)}$ with $i'\geq1$, and set $\tp_i' = P_i \tp_i$ for $i=0,1$, $\td' = \tp_0' \cdot \tp_1'$,  and $m =P_0P_1$.
 We define a unital embedding $\psi_n: C(\tX_n) \otimes M_{\td} \rightarrow C(\tX_{n+1}) \otimes M_{\td'} $ by dividing into the following three cases depending on $x\in \tX_{n+1}$: For $f\in C(\tX_n) \otimes M_{\td}$ and
for $x=\iota_{0,i}^{(n+1)} (t) \in I_{0,i}^{(n+1)}$, $i=0,1$, and $t\in [0,2\pi]$ we define
\[ \psi_n (f) (x) = f(\iota_{0,i}^{(n)} (t)) \otimes 1_{m};\]
For $x = \tiota_{i,j'}^{(n+1)} ((-1)^it)\in \tI_{i,j'}^{(n+1)}$, $i=0,1$, $j'=1,2,...,r_0^{(n)}$, and $t\in [0,\pi]$ we define
\[\psi_n  (f) (x) = \bigoplus _{j=1}^{r_0^{(n)}} f(\tiota_{\natural(i,j,j'),j}^{(n)}((-1)^{\natural(i,j,j')}t)) \otimes 1_{|\mu_{j',j}|} \oplus f(c_n)\otimes  1_{m_{0,j'}}, \]
where $\natural(i,j,j') \in \{0,1\}$ is defined by $\arg(\mu_{j',j})(-1)^i = (-1)^{\natural(i,j,j')}$  and where $m_{0,j'}=m- \sum_{j=1}^{r_0^{(n)}} |\mu_{j',j}|$; For $x= \iota_{i',j'}^{(n+1)}(t) \in I_{i',j'}^{(n+1)}$, $(i',j')\in S^{(n+1)}$ with $i'\geq 1$ and $t \in [0,2\pi] $, we define
\begin{eqnarray}
\psi_n(f) (x) &=& \bigoplus _{j=1}^{r_0^{(n)}} (f(\tiota_{0,j}^{(n)}(t/2)) \otimes 1_{\nu_{i',j',0,j}} \oplus f(\tiota_{1,j}^{(n)}(-t/2))\otimes 1_{d_{i',j'}^{(n+1)} - \nu_{i',j',0,j}})
\nonumber \\
&\oplus& \bigoplus_{j=1}^{r_{i'}^{(n)}} (f (\iota_{i'j}^{(n)} (t)) \otimes 1_{\nu_{j',i',j}} \oplus f(\iota_{0,\natural(i',j')}^{(n)} (t)) \otimes 1_{d_{i',j'}^{(n+1)}-\nu_{j',i',j}}) \nonumber \\
&\oplus& f(c_n) \otimes 1_{m-(r_0^{(n)}+ r_{i'}^{(n)})d_{i',j'}^{(n+1)}}. \nonumber 
\end{eqnarray}
Note that when $p_{i'} \nmid |G_{n,1}|$ we have defined $ r_{i'}^{(n)}=0$.
It is not so hard to show that $\psi_n (f) (x) \rightarrow f(c_n) \otimes 1_{m}$, $x\rightarrow c_{n+1}$ for $f \in C(\tX_n)\otimes M_{\td} $. Thus it follows that $\psi_n (C(\tX_n)\otimes M_{\td})\subset C(\tX_{n+1})\otimes M_{\td'}$.  

Recall that $A_{n} = A(G_{n},\tp_0, \tp_1)$ was defined in terms of automorphisms $\epsilon_{j}^{(n)}$ and also depend on unital embeddings $\epsilon_{i,j}^{(n)}$. We shall show that there exist automorphisms $\epsilon_{j'}^{(n+1)}$ of $M_{\td'}$ and unital embeddings $\epsilon_{i',j'}^{(n+1)}$ which satisfy the condition in (0,0). Let $\epsilon_{j'}^{(n+1)}$ be an extension of $(\bigoplus_{j=1}^{r_0^{(n)}} (\epsilon_j^{(n)})^{\arg (\mu_{j',j})} \otimes 1_{|\mu_{j',j}|})\oplus\id_{M_{\td}}\otimes 1_{m_{0,j'}}$, then we have $\epsilon_{j'}^{(n+1)}(\psi_n(f)(b_{1,j'}^{(n+1)}))=\psi_n(f)(b_{0,j'}^{(n+1)})$ for $j'=1,2,...,r_0^{(n+1)}$. 
Since $\tp_{\natural (i)}' | m \tp_{\natural(i)} $, $i=0,1$, $f(a_{0,i}^{(n)}) \in \epsilon_{0,i}^{(n)}(M_{\tp_i})$ for any $f \in A_n$, and $\psi_n (f) (a_{0,i}^{(n+1)}) = f(a_{0,i}^{(n)}) \otimes 1_m$ there exists a unital embedding $\epsilon_{0,i}^{(n+1)}:M_{\tp_{i}'}\hookrightarrow M_{\td'}$ such that 
\[ \psi_n (f) (a_{0,i}^{(n+1)}) \in \epsilon_{0,i}^{(n+1)} (M_{\tp_{i}'}), \quad f \in A_n,\ i=0,1.\]
By the definition of $\psi_n$, it follows that for $f\in A_n,$ 
\begin{eqnarray} 
\psi_n (f) (a_{i',j'}^{(n+1)}) &=& \bigoplus_{j=1}^{r_0^{(n)}} \epsilon_j^{(n)} (f(b_{1,j}^{(n)}))\otimes 1_{\nu_{i',j',0,j}}\oplus f(b_{1,j}^{(n)})\otimes 1_{d_{i',j'}^{(n+1)}-\nu_{i',j',0,j}} \nonumber \\
&\oplus& \bigoplus_{j=1}^{r_{i'}^{(n)}} f(a_{i',j}^{(n)})\otimes 1_{\nu_{j',i',j}} \oplus f(a_{0,\natural(i',j')}^{(n)}) \otimes 1_{d_{i',j'}^{(n+1)} -\nu_{j',i',j}} \nonumber \\  
&\oplus& f(c_n) \otimes 1_{m-(r_0^{(n)} +r_{i'}^{(n)})d_{i',j'}^{(n+1)}},\nonumber
\end{eqnarray} 
Since $d_{i',j'}^{(n+1)} | d_{i',j}^{(n)} \nu_{j',i',j}$ for any $j=1,2,...,r_{i'}^{(n)}$ and since $f(a_{i',j}^{(n)}) \in \epsilon_{i',j}^{(n)}(M_{\td / d_{i',j}^{(n)}})$ for  $f \in A_n$, we obtain a unitary $u_{i',j',j}$ in $M_{\td \nu_{j',i',j}}$ such that 
\[\Ad u_{i',j',j} (f(a_{i',j}^{(n)})\otimes 1_{\nu_{j',i',j}}) \in M_{\td\nu_{j',i',j}/d_{i',j'}^{(n+1)}}\otimes 1_{d_{i',j'}^{(n+1)}},\quad f\in A_n,\ p_{i'} |\ |G_{n,1}|. \]
When $p_{i'} \nmid |G_{n,1}| $ we have defined $\nu_{j',i',j}=0$.  
When $p_{i'} \mid |G_{n,1}|$ we have $d_{i',j}^{(n)} \mid \tp_{\natural(\natural(i',j))}$, then it follows  that $d_{i',j'}^{(n+1)} \mid \tp_{\natural(\natural(i',j))} (d_{i',j'}^{(n+1)} - \nu_{j',i',j}),$ and from $f(a_{0,\natural(i',j)}^{(n)}) \in \epsilon_{0,\natural(i',j)}^{(n)}(M_{\tp_{\natural(i',j)}})$ for $f \in A_n$ we obtain a unitary $v_{i',j',j}$ in $M_{\td(d_{i',j'}^{(n+1)}- \nu_{j',i',j})}$ such that 
\[\Ad v_{i',j',j} (f(a_{0,\natural(i',j')}^{(n)})\otimes 1_{d_{i',j'}^{(n+1)}-\nu_{j',i',j}}) \in M_N \otimes 1_{d_{i',j'}^{(n+1)}}, \quad f \in A_n, \]
where $N=\td(d_{i',j'}^{(n+1)} - \nu_{j',i',j})/d_{i',j'}^{(n+1)}$.  Hence there exists a unital embedding $\epsilon_{i',j'}^{(n+1)}:M_{\td'/d_{i',j'}^{(n+1)}}\hookrightarrow M_{\td'}$ such that 
\[ \psi_n(f) (a_{i',j'}^{(n+1)}) \in \epsilon_{i',j'}^{(n+1)} (M_{\td'/d_{i',j'}^{(n+1)}}), \quad f \in A_n.\]
Whence we conclude that $\psi_n (A_n) \subset A(G_{n+1},\tp_0',\tp_1')$.

{\noindent{\bf Proof\ of\ ($\psi_n|_{A_n})_* \cong \gamma_n$}. }
 In Proposition \ref{Prop1}, we have identified $[f_{0,j}^{(n)}]_{K_1(A_n)}$ and $[f_{i,j}^{(n)}]_{K_1(A_n)}$ with $g_{0,j}^{(n)}$ and $ g_{i,j}^{(n)}.$ Set $A_{n+1}=A(G_{n+1},\tp_0',\tp_1')$ and
\begin{eqnarray}
f_{0,j}'(x)=\left\{ \begin{array}{ll}
\exp ( 2it e_{0,0} ) ,\quad & x= \tiota_{0,j}^{(n)} (t) \in \tI_{0,j}^{(n)},\ 0\leq t \leq \pi , \\
1_{\td}, \quad  & {\rm otherwise}, 
\end{array} \right.\nonumber
\end{eqnarray}
where $e_{0,0}$ is a minimal projection of $M_{\td}$. Remark that $f_{0,j}' \in U(A_n)$ and $[f_{0,j}']_{K_1(A_n)} = [f_{0,j}^{(n)}]=g_{0,j}^{(n)}$. By the definition of $\psi_n$ it follows that $\psi_n(f_{0,j}')(x)$
\begin{eqnarray}
=\left\{ \begin{array}{lll}
\exp ( 2it e_{0,0}\otimes 1_{|\mu_{j',j}|}) ,\quad & x= \tiota_{\natural(0,j,j'),j'}^{(n+1)} (\arg(\mu_{j',j})t),\ 0\leq t \leq \pi, \\
\exp (ite_{0,0} \otimes 1_{\nu_{i',j',0,j}}), \quad & x = \iota_{i',j'}^{(n+1)}(t) \in I_{i',j'}^{(n+1)},\ i' \geq 1,\ 0\leq t \leq 2\pi, \\ 
1_{\td'}, \quad  & {\rm oterwise},  
\end{array} \right.\nonumber
\end{eqnarray}
which implies that $\psi_n (f_{0,j}') \in U(A_{n+1})$ and  \[ [ \psi_n (f_{0,j}')]_{K_1(A_{n+1})} = \sum _{j'} \mu_{j',j} g_{0,j'}^{(n+1)} + \sum_{i'}\sum_{j'} \nu_{i',j',0,j} g_{i',j'}^{(n+1)}.\]
For $(i,j)\in S^{(n)}$ with $i\geq 1 $ it follows that
\begin{eqnarray}
\psi_n(f_{i,j}^{(n)})(x)=\left\{ \begin{array}{ll}
\exp (ite_{0,0} \otimes 1_{\nu_{j',i,j}}), \quad & x = \iota_{i,j'}^{(n+1)}(t) \in I_{i,j'}^{(n+1)},  \\
1_{\td'}, \quad  & {\rm otherwise} ,
\end{array} \right.\nonumber
\end{eqnarray} 
which implies that $\psi_n (f_{i,j}^{(n)}) \in U(A_{n+1})$ and  
\[ [ \psi_n (f_{i,j}^{(n)})]_{K_1(A_{n+1})} = \sum_{j'} \nu_{j',i,j} g_{i,j'}^{(n+1)}.\]
Hence we have that $( \psi_n)_* \cong \gamma_n$ .

{\noindent{\bf Proof\ of\ (0,1)}. }
Let $\hpsi_{n,\delta}: \tX_{n+1} \rightarrow \tX_n,$ $\delta=1,2,...,m$ be the continuous maps such that $\displaystyle \bigoplus_{\delta=1}^{m} f\circ \hpsi_{n,\delta} (x) = \psi_n (f) (x),$  $f \in C(\tX_n)\otimes M_{\td}$. 
 Since $\gamma_n$ is injective, for $j_0 \in \{1,2,...,r_0^{(n)}\}$, there exists  $j_0' \in \{1,2,...,r_0^{(n+1)}\}$ such that $\mu_{j_0',j_0} \neq 0$. Set $\displaystyle \delta_0=\sum_{j=1}^{j_0-1}|\mu_{j_0',j}|+1$. 
By the definition of $\psi_n$ at $x=\tiota_{i,j_0'}((-1)^it)$, we have that
$\hpsi_{n,\delta_0}(\tI_{0,j_0'}^{(n+1)} \cup \tI_{1,j_0'}^{(n+1)})= \tI_{0,j_0}^{(n)} \cup \tI_{1,j_0}^{(n)}$ and that if $\psi_n(f)^{-1}(\{ 0 \}) \cap (\tI_{0,j_0'}^{(n+1)} \cup \tI_{1,j_0'}^{(n+1)})$ is $\varepsilon$-dense then $\hpsi_{n,\delta_0}(\psi_n(f)^{-1} ( \{ 0 \}) \cap (\tI_{0,j_0'}^{(n+1)} \cup \tI_{1,j_0'}^{(n+1)}))$ is $\varepsilon$-dense in $\tI_{0,j_0}^{(n)} \cup \tI_{1,j_0}^{(n)}$. Since $f^{-1} (\{ 0 \}) \cap (\tI_{0,j_0}^{(n)} \cup \tI_{1,j_0}^{(n)}) \supset \hpsi_{n,\delta_0} (\psi_n(f)^{-1} (\{ 0 \}) \cap (\tI_{0,j_0'}^{(n+1)} \cup \tI_{1,j_0'}^{(n+1)}))$, it follows that if $\psi_n(f)^{-1}(\{ 0 \}) \cap (\tI_{0,j_0'}^{(n+1)} \cup \tI_{1,j_0'}^{(n+1)})$ is $\varepsilon$-dense then $f^{-1}(\{ 0 \}) \cap (\tI_{0,j_0}^{(n)} \cup \tI_{1,j_0}^{(n)}) $ is $\varepsilon$-dense. 
Similarly, since $\gamma_n$ is injective, for $(i_0,j_0) \in S^{(n)}$ there exists $(i_0, j_0') \in S^{(n+1)}$ such that $\nu_{j_0',i_0,j_0} \neq 0 .$ Set $\displaystyle \delta_0= (r_0^{(n)}+ j_0-1) d_{i_0,j_0'}^{(n+1)} +1.$ By the defitnition of  $\psi_n$ at $x=\iota_{i_0,j_0'}^{(n+1)} (t) $, we have that $\hpsi_{n,\delta_0}( I_{i_0,j_0'}^{(n+1)}) = I_{i_0, j_0}^{(n)}$ and that if $\psi_n (f)^{-1} (\{ 0 \}) \cap I_{i_0,j_0'}^{(n+1)}$ is $\varepsilon$-dense then $\hpsi_{n,\delta_0}(\psi_n(f)^{-1} (\{ 0 \}) \cap I_{i_0, j_0'}^{(n+1)}) \subset f^{-1} (\{ 0 \}) \cap I_{i_0,j_0}^{(n)}$ is $\varepsilon$-dense.
Considering a neighborhood of $c_n$, we notice that if $\psi_n (f)^{-1} ( \{ 0 \}) \cap (\tX_{n+1} \setminus (I_{0,0}^{(n+1)} \cup I_{0,1}^{(n+1)})) \neq \phi$
then $f(c_n) = 0$. If $\psi_n(f)^{-1} ( \{ 0 \}) \cap (\tX_{n+1}\setminus (I_{0,0}^{(n+1)} \cup I_{0,1}^{(n+1)})) = \phi$ and $\psi_n(f)^{-1} (\{ 0\})$  is $\varepsilon$-dense then there exists $(0,j) \in S^{(n+1)}$ such that $d(\psi_n(f)^{-1}(\{ 0 \}) \cap I_{0,j}^{(n+1)}, c_{n+1}) \leq \varepsilon ,$ which induces that $d(f^{-1}(\{ 0 \})  \cap I_{0,j}^{(n)}, c_n) \leq \varepsilon,$ $j=0$ or $1$. Hence, combining the above arguments, we have (0,1).

{\noindent{\bf Proof\ of\ (0,2)}. }
Let $\hpsi_{n,\delta}: \tX_{n+1} \rightarrow \tX_n,$ $\delta=1,2,...,m$ be the continuous maps such that $\displaystyle \bigoplus_{\delta=1}^{m} f\circ \hpsi_{n,\delta} (x) = \psi_n (f) (x),$  $f \in C(\tX_n)\otimes M_{\td}$. 
By the definition of $\psi_n$, we have that 
\[ d (\hpsi_{n,\delta} (x), \hpsi_{n,\delta}(y) ) \leq d(x,y), \]
for any $x,y\in \tX_{n+1}$ and $\delta=1,2,...,m$. Then for $x$ and $y\in \tX_{n+1}$ with $d(x,y) \leq \varepsilon$ and $f \in C(\tX_n) \otimes M_{\td}$, it follows that 
\begin{eqnarray}
& &|\Tr_{\td'} (\psi_n(f) (x) - \psi_n(f) (y) ) |  \nonumber \\
& &\leq \frac{1}{m} \sum_{\delta=1}^{m} |\Tr_{\td} (f(\hpsi_{n,\delta} (x))- f(\hpsi_{n,\delta}(y))) | \nonumber \\
& &\leq \max \{ |\Tr(f(x) -f(y))| ;\ x,y \in \tX_n ,\ d(x,y) \leq \varepsilon \}.\nonumber\end{eqnarray}
Thus we have (0,2). 

{\noindent{\bf Proof\ of\ (0,3)}. }
By the definition of $\psi_n$, for $h^{(n)}_k \in C(\tX_n)_+^1,$ $n,$ $k \in \N$ with $h^{(n)}_k \searrow \chi(E_n)$ $(k\rightarrow \infty)$, $\psi_n$ satisfies  that 
\[ \lim_{k \rightarrow \infty} \psi_n(h_k^{(n)} \otimes 1_{\td}) \leq \chi( E_{n+1})\otimes 1_{\td'}.\] Then we have (0,3). This completes the proof of the proposition.
\eprf

\section{ Simple projectionless $C^*$-algebras with a unique tracial state}\label{Sec3}
In this section, we construct the inductive limit $C^*$-algebra of a sequence of dimension drop algebras that is simple and has a unique tracial state.

For the simplicity of  $C^*$-algebras and the uniqueness of tracial states, we prepare the following proposition. 

\begin{proposition}\label{Prop3}
Let $G_n$ be a finitely generated abelian group, $\tp_0$ and $\tp_1$ relatively prime numbers such that $d_{i,j}^{(n)} | \td= \tp_0 \cdot \tp_1$ for any $(i,j)\in S^{(n)}$ with $i \geq 1$, $A_n$ a dimension drop algebra for $G_n$, $\tp_0$, and $\tp_1$.
Then there exist relatively prime numbers $\tp_0'$ and $\tp_1' $ such that $\tp_0| \tp_0'$ and $\tp_1 | \tp_1'$, automorphisms $\{ \epsilon_j^{(n,0)} \in Aut (M_{\td'});$ $j=1,2,...,r_0^{(n)}\}$, where $\td'= \tp_0'\tp_1'$, unital embeddings $\{\epsilon_{i,j}^{(n,0)}:M_{\td'/d_{i,j}^{(n)}} \hookrightarrow M_{\td'};$ $(i,j)\in S^{(n)}\}$, where $d_{0,i}^{(n)} = \tp_{\natural(i)}'$ for $i=0,1$, and a unital embedding $\varphi_{n,0} : C(\tX_n)\otimes M_{\td} \hookrightarrow C(\tX_{n})\otimes M_{\td'}$ satisfying the following conditions:
\begin{description}\label{prop:embedding}
\item[(0)]\label{(0)} $\varphi_{n,0} (A_n) \subset A(G_n,\tp_0',\tp_1')$ and $(\varphi_{n,0}|_{A_n})_* \cong \id_{G_n}$ under the identification which is defined in Proposition \ref{Prop1}, where $A(G_n,\tp_0',\tp_1')$ is defined by $\epsilon_j^{(n,0)}$ and $\epsilon_{i,j}^{(n,0)}$,
\item[(1)]\label{(1)} if $f\in C(\tX_n)\otimes M_{\td}$ satisfies that $\varphi_{n,0}(f)^{-1} (\{ 0 \}) $ is $\varepsilon$-dense in $\tX_n$ then $f^{-1} (\{0 \}) $ is $1/2 \varepsilon$-dense in $\tX_n$,
\item[(2)]\label{(2)} for any $f\in C(\tX_n)\otimes M_{\td}$ and $\varepsilon >0 $, 
\[ \max \{ |\Tr_{\td'}(\varphi_{n,0}(f) (x)- \varphi_{n,0} (f)(y))| ;\ x,y  \in \tX_n,\ d(x,y) \leq \varepsilon \}\] 
\[ \leq \max \{ |\Tr_{\td} (f(x)-f(y))| ;\ x,y \in \tX_n,\ d(x,y) \leq 1/2 \varepsilon\},\] 
where $d$ is the canonical metric on $\tX_n$,
\item[(3)]\label{(3)} 
$\Tr_{\td'} ( \varphi_{n,0} (\chi (E_n) \otimes 1_{\td}) (x)) \leq 1/2 \Tr (\chi(E_n) \otimes 1_{\td'}(x)), \quad x \in \tX_n.$
\end{description} 
\end{proposition}
\bprf
Let $G_n$,  $\tp_0$, $\tp_1$, and $\td$ be as in Section \ref{Sec2} (i.e., by
$r_0^{(n)} \in \Z_+ ,$ $r_i^{(n)},$ $k_{i,j}^{(n)} \in \N $ for $i \geq 1 $ with $p_i |\ |G_{n,1}|$, and $d_{i,j}^{(n)} = p_i^{k_{i,j}^{(n)}}$, we have represented $G_n$
as   
\[G_n \cong \Z^{r_0^{(n)}} \oplus \bigoplus_{p_i |\ |G_{n,1}|}\ \bigoplus_{j=1}^{r_i^{(n)}} \Z_{d_{i,j}^{(n)}},\] 
$(\tp_0, \tp_1) =1$, $\td= \tp_0 \cdot \tp_1$, and $ d_{i,j}^{(n)} | \td$ for any $(i,j)\in S^{(n)}$ with $i \geq 1$),  
and let $A_n$ be a dimension drop algebra for $G_n$, $\tp_0$, and $\tp_1$.

First, we define $\tp_i'$, $i=0,1$ as follows. Let  
$P_i \in \N$, $i=0,1$ be such that $(P_0 \tp_0 , P_1 \tp_1 )=1$, $P_i > \td$, and 
\[P_0 \equiv   P_1 \equiv 1 \pmod{\td}\]
( $P_0=a'\td +1,$ $P_1=b'P_0\td +1,$ $a',$ $b'\in\N$).
We define $\tp_i' = P_i \tp_i $ for $i=0,1$, and $m= P_0P_1$.
By the definition of $P_0$ and $P_1$, we obtain $a$, $b$ and $c \in \N$ satisfying that
\[ m =a \tp_0' + P_0=  b\tp_1' + P_1 =  c\td +1.\]
Set 
\[ \kappa = P_0 + P_1 + \sum_{p_i|\ |G_{n,1}|} r_i^{(n)} ,\]
then, taking large $P_i$, we may assume that $ m \geq 2( 3 r_0^{(n)}
+\kappa)$.

We shall define $\varphi_{n,0}$ as the direct sum of $r_0^{(n)}$ homomorphisms of $C(\tX_n)\otimes M_{\td}$ into $C(\tX_n)\otimes M_{\td} \otimes M_3$, $\kappa$ homomorphisms of $C(\tX_n)\otimes M_{\td} $ into $C(\tX_n)\otimes M_{\td}$, and the homomorphism $\xi_{c}$ of $C(\tX_n)\otimes M_{\td}$ into $C(\tX_n)\otimes M_{\td}\otimes M_{m-(3r_0^{(n)}+\kappa)}$ defined by $f\mapsto f(c_n)\otimes 1_{m-(3r_0^{(n)}+\kappa)}$. The first type of homomorphisms, denoted by $\eta_j$, will be constructed to preserve the free part of $G_n$ and  the second type denoted by $\xi_{\delta}$ will be constructed to preserve the torsion subgroup of $G_n$. We define a $*$-homomorphism $\eta: C(\tX_n) \otimes M_{\td} \rightarrow C(\tX_n) \otimes M_{3r_0^{(n)}\td}$ as follows. 
Set 
\begin{eqnarray}
D (t) &=& \diag ( z(\frac{t}{3}),z(\frac{t+2\pi}{3}),z(\frac{t-2\pi}{3})) \in U(M_{3}),\quad t \in[ -\pi, \pi],\nonumber \\
  F(t) &=& \diag( 0,  \frac{2\pi -t}{3}, \frac{t-2\pi}{3})\in M_3,\quad t\in[0,2\pi]. \nonumber 
\end{eqnarray}
and let $V \in U(M_{3}) $ be such that
\[ V D(-\pi) V^* = D(\pi). \]
For $f \in C(\tX_n)\otimes M_{\td}$ and $j \in \{1,2,...,r_0^{(n)}\}$ let  $f \circ \tau_j^{(n)} (D (t))$
\begin{eqnarray}
=\left\{ \begin{array}{ll}
 \diag (f\circ \tiota_{0,j} (\frac{t}{3}), f\circ \tiota_{0,j} (\frac{t+ 2\pi}{3}),f\circ \tiota_{1,j} (\frac{t-2\pi}{3})),\quad  t \in [0, \pi]  \\
 \diag (f\circ \tiota_{1,j} (\frac{t}{3}), f\circ \tiota_{0,j} (\frac{t+ 2\pi}{3}), f\circ \tiota_{1,j} (\frac{t- 2\pi}{3})),\quad  t\in [-\pi, 0], \\ 
\end{array} \right.\nonumber
\end{eqnarray}
 and let $f \circ \tau_j^{(n)}\circ z (F (t))$
\[ = \diag (f(c_n) , f\circ \tiota_{0,j} (\frac{2\pi-t}{3}), f\circ \tiota_{1,j}(\frac{t-2\pi}{3})), \quad t\in [0, 2\pi].\]
Remark that for $f \in  A_n $,
\[ \Ad 1_{\td} \otimes V\cdot (1_{\td}\oplus 1_{\td}\oplus v_j^{(n)}) (f\circ \tau_{j}^{(n)} (D(-\pi))) = f\circ\tau_j^{(n)}(D(\pi)),\]
where $v_j^{(n)}$ is a unitary in $M_{\td}$ satisfying that $\Ad v_j^{(n)} (f(b_{1,j}^{(n)}))= f(b_{0,j}^{(n)})$ for any $f\in A_n$.
Set $v_j' = 1_{\td}\otimes V \cdot (1\oplus 1\oplus v_j^{(n)})$.

Define $*$-homomorphisms $\eta_j : C(\tX_n) \otimes M_{\td} \rightarrow C(\tX_n) \otimes M_{\td}\otimes M_{3}$, $j=1,2,..., r_0^{(n)} $ by  
\begin{eqnarray}
\eta_{j}(f)(x)=\left\{ \begin{array}{lll}
 f\circ\tau_{j}^{(n)} (D ((-1)^it)), \quad & x= \tiota_{i,j}^{(n)}((-1)^it),\ t\in[0,\pi],\ i=0,1,  \\
 f\circ\tau_{j}^{(n)} (D (0)), \quad & x \in \tI_{0,\nu}^{(n)} \vee \tI_{1,\nu}^{(n)},\ \nu \neq j, \\
f\circ \tau_{j}^{(n)}\circ z (F (t)), \quad & x = \iota_{\mu,\nu}^{(n)} (t) \in I_{\mu,\nu}^{(n)},\ (\mu,\nu)\in S^{(n)}, 
\end{array} \right.\nonumber
\end{eqnarray}
for $f\in C(\tX_n)\otimes M_{\td}$.
Remark that \[ \eta_j(f) (a_{\mu,\nu}^{(n)}) = f (c_n) \otimes 1_{3},\quad f \in C(\tX_n) \otimes M_{\td},\ (\mu,\nu) \in S^{(n)}, \]
and that \[ \Ad v_j' ( \eta_j (f) ( b_{1,j}^{(n)}) )= \eta_j (f) (b_{0,j}^{(n)}),\quad  f\in A_n,\ j=1,2,...,r_0^{(n)}.\]
We define a $*$-homomorphism  $\eta: C(\tX_n) \otimes M_{\td} \rightarrow C(\tX_n) \otimes M_{3r_0^{(n)}\td}$ by $\eta = \bigoplus_{j=1}^{r_0^{(n)}} \eta_j$ .

We define a $*$-homomorphism $\xi_{\delta}^{(n)}: C(\tX_n) \otimes M_{\td} \rightarrow C(\tX_n) \otimes M_{\td},$ $\delta =1,2,..., \kappa$ by the following cases :
when $f\in C(\tX_n) \otimes M_{\td} $ and $1\leq \delta \leq P_1 $
\begin{eqnarray}
\xi_{\delta}^{(n)}(f)(x)=\left\{ \begin{array}{lll}
f(\iota_{0,0}^{(n)}(\pi + \frac{t}{2})), \quad & x = \iota_{0,0}^{(n)} (t) \in I_{0,0}^{(n)}, \\
f(\iota_{0,0}^{(n)}(\pi - \frac{t}{2})), \quad & x= \iota_{\mu,\nu }^{(n)} (t) \in I_{\mu,\nu}^{(n)}, \ (\mu,\nu )\neq (0,0), \\
f(\iota_{0,0}^{(n)}(\pi)), \quad & x\in \tI_{0,j}^{(n)}\cup \tI_{1,j}^{(n)},\ j=1,2,...,r_0^{(n)},
\end{array} \right.\nonumber
\end{eqnarray} 
when $P_1+1 \leq\delta \leq P_0+P_1 $
\begin{eqnarray}
\xi_{\delta}^{(n)}(f)(x)=\left\{ \begin{array}{lll}
f(\iota_{0,1}^{(n)}(\pi + \frac{t}{2})), \quad & x = \iota_{0,1}^{(n)} (t) \in I_{0,1}^{(n)}, \\
f(\iota_{0,1}^{(n)}(\pi - \frac{t}{2})), \quad & x= \iota_{\mu, \nu}^{(n)} (t) \in I_{\mu, \nu}^{(n)}, \ (\mu, \nu)\neq (0,1), \\
f(\iota_{0,1}^{(n)} (\pi)), \quad & x\in \tI_{0,j}^{(n)}\cup \tI_{1,j}^{(n)},\ j=1,2,...,r_0^{(n)},
\end{array} \right.\nonumber
\end{eqnarray} 
when $\displaystyle \delta = P_0+ P_1 + \sum_{\mu< i,\ p_{\mu}|\ |G_{n,1}|} r_{\mu}^{(n)} +j $, $p_i |\ |G_{n,1}|$, and $j=1,2,...,r_i^{(n)}$   
\begin{eqnarray}
\xi_{\delta}^{(n)} (f) (x)=\left\{ \begin{array}{lll}
f(\iota_{i,j}^{(n)}(\pi + \frac{t}{2})), \quad & x = \iota_{i,j}^{(n)} (t) \in I_{i,j}^{(n)}, \\
f(\iota_{i,j}^{(n)}(\pi - \frac{t}{2})), \quad & x= \iota_{\mu, \nu }^{(n)} (t) \in I_{\mu,\nu}^{(n)}, \ (\mu , \nu)\neq (i,j), \\
f(\iota_{i,j}^{(n)} (\pi)), \quad & x\in \tI_{0,j}^{(n)}\cup \tI_{1,j}^{(n)},\ j=1,2,...,r_0^{(n)},
\end{array} \right.\nonumber
\end{eqnarray} 
We define a $*$-homomorphism $\xi: C(\tX_n) \otimes M_{\td} \rightarrow C(\tX_n) \otimes M_{\kappa \td}$ by $\xi = \displaystyle \bigoplus_{\delta =1}^{\kappa} \xi_{\delta} $. 

Set $\td' = \tp_0'\tp_1' (= m \td )$, and let $A_n$ be defined by automorphisms $\epsilon_{j}^{(n)}(=\Ad v_j^{(n)})$ and unital embeddings $\epsilon_{i,j}^{(n)}$. Define a $*$-homomorphism $\xi_c: C(\tX_n) \otimes M_{\td} \rightarrow C(\tX_n) \otimes
M_{(m-(3r_0^{(n)}+\kappa))\td}$ by $\xi_c (f) (x) = f(c_n) \otimes 1_{m-(3r_0^{(n)} +\kappa)}$, $ x \in \tX_n$, and define a unital $*$-homomorphisms $\varphi_{n,0}: C(\tX_n)\otimes M_{\td} \hookrightarrow C(\tX_n)\otimes M_{\td'} $ by 
\[ \varphi_{n,0} = \eta \oplus \xi \oplus \xi_c.\]
We shall show that there exist automorphisms $\epsilon_j^{(n,0)}$ and unital embeddings $\epsilon_{i,j}^{(n,0)}$ satisfying the condition in (0). Set $v_j=1_{(j-1)\td} \oplus v_j' \oplus 1_{(m-j)\td} \in U (M_{\td'})$ and $\epsilon_j^{(n,0)} = \Ad v_j$, then it follows  that 
\[ \epsilon_j^{(n,0)} (\varphi_{n,0} (f) (b_{1,j}^{(n)})) = \varphi_{n,0} (f) (b_{0,j}^{(n)}), \quad f \in A_n,\ j=1,2,...,r_0^{(n)}.\] 
From the remark of $\eta$ and the definition of $\xi$ it follows that 
\[\varphi_{n,0} (f) (a_{0,i}^{(n)} ) = f(c_n) \oplus \cdots \oplus f(c_n) \oplus f(a_{0,i}^{(n)}) \otimes 1_{P_{\natural(i)}} \oplus f(c_n) \oplus \cdots 
\oplus f(c_n),\]
for $i=0,1$. Since $f(a_{0,i}^{(n)}) \in \epsilon_{0,i}^{(n)}(M_{\tp_i})$ for $f \in A_n$, $\tp_{\natural(i)}' = P_{\natural(i)}\tp_{\natural(i)} $, and $\tp_{\natural(i)}' | m -P_{\natural(i)}$ there exist unital embeddings  $\epsilon_{0,i}^{(n,0)} : M_{\tp_i'} \hookrightarrow M_{\td'}$, $i=0,1,$ such that 
\[\varphi_{n,0} (f) (a_{0,i}^{(n)})  \in \epsilon_{0,i}^{(n,0)}(M_{\tp_i'}),\ {\rm for\ any\ }f \in A_n.\]
Similarly, it follows that 
\[ \varphi_{n,0} (f) (a_{i,j}^{(n)}) = f(c_n) \oplus \cdots \oplus f(c_n) \oplus f(a_{i,j}^{(n)})\oplus f(c_n) \oplus \cdots \oplus f(c_n),\]
for any $(i,j)\in S^{(n)}$ with $i\geq 1$.
Since $f(a_{i,j}^{(n)}) \in \epsilon_{i,j}^{(n)} (M_{\td / d_{i,j}^{(n)}})$ and $d_{i,j}^{(n)}|\ \td\ |\ m - 1$, there exist unital embeddings $\epsilon_{i,j}^{(n,0)}: M_{\td'/ d_{i,j}^{(n)}}\hookrightarrow M_{\td'}$ such that  
\[ \varphi_{n,0}(f) (a_{i,j}^{(n)}) \in \epsilon_{i,j}^{(n,0)}(M_{\td'/ d_{i,j}^{(n)}}),\ {\rm for \ any\ } f \in A_n .\]
Hence  we conclude that $\varphi_{n,0} (A_n) \subset A(G_n,\tp_0',\tp_1') $.

\

 {\noindent{\bf Proof\ of\ $(\varphi_{n,0}|_{A_{n}})_*=\id_{G_n}$}.\ }
Define $f_{0,j}' \in U(A_n)$ for $j=1,2,...,r_0^{(n)}$ by
\begin{eqnarray}
f_{0,j}' (x)=\left\{ \begin{array}{ll}
\exp (6it e_{0,0}) ,\quad & x= \tiota_{0,j}^{(n)}(t),\ 0 \leq t \leq \pi/3, \nonumber \\
1_{\td}, \quad  & {\rm otherwise}, \\
\end{array} \right.\nonumber
\end{eqnarray} 
and  $f_{i,j}' \in U(A_n) $ for  $(i,j) \in S^{(n)} $ with $i\geq 1$ by
\begin{eqnarray}
f_{i,j}' (x)=\left\{ \begin{array}{ll}
\exp (2i(t-\pi) e_{0,0}) ,\quad & x= \iota_{i,j}^{(n)}(t),\ \pi \leq t \leq 2\pi, \nonumber \\
1_{\td}, \quad  & {\rm otherwise}. \\
\end{array} \right.\nonumber
\end{eqnarray} 
Remark that $[f_{0,j}']_{K_1(A_n)} = [ f_{0,j}^{(n)} ] = g_{0,j}^{(n)},$ $j=1,2,...,r_0^{(n)}$, and $[ f_{i,j}']_{K_1(A_n)} = [ f_{i,j}^{(n)}]  = g_{i,j}^{(n)}$, $(i,j) \in S^{(n)}$ with $i \geq 1$.  
By the definitions of $\eta,$ $\xi,$ and $\xi_c$, it follows that $\eta_{\nu}(f_{0,j}') \equiv 1_{3\td}$, $\nu \neq j$, $(\xi \oplus \xi_c)(f_{0,j}') \equiv 1_{(m-3r_0^{(n)})\td}$, and 
\begin{eqnarray}
\eta_j(f_{0,j}' )(x)=\left\{ \begin{array}{lll}
\exp (2 i t e_{0,0}) ,\quad & x= \tiota_{0,j}^{(n)}(t),\ 0\leq t \leq \pi, \nonumber \\
\exp(-2 i te_{0,0}), \quad & x= \iota_{i,j}^{(n)} (t), \ \pi \leq t \leq 2\pi, \nonumber \\
1_{3\td}, \quad  & {\rm otherwise}. \\
\end{array} \right.\nonumber
\end{eqnarray} 
Then we have that 
\[ [\varphi_{n,0} (f_{0,j}')]_{K_1 (A(G_n,\tp_0', \tp_1'))} = [ ( \eta \oplus \xi \oplus \xi_c )(f_{0,j}')]= g_{0,j}^{(n)}.\]
By the definition of $\eta$ it follows that $ \eta_{\nu} (f_{i,j}') \equiv 1_{3\td}, \quad \nu =1,2,..., r_0^{(n)},$ and by the definition of $\xi _{\delta}^{(n)}$ it follows that $ \xi_{\delta}^{(n)} (f_{i,j}') \equiv  1_{\td}$ for $\delta \neq P_0 + P_1 + \sum_{\mu<i, \ p_{\mu}|\ |G_{n,1}|} r_\mu^{(n)} + j$ and $\xi_{\delta}(f_{i,j}') = f^{(n)}_{i,j}$ when $\delta = P_0 + P_1 + \sum r_\mu^{(n)} + j$. Then we have that 
\[ [\varphi_{n,0} (f_{i,j}')]_{K_1 (A(G_n,\tp_0',\tp_1'))} = [( \eta \oplus \xi \oplus \xi_c) (f_{i,j}') ] =g_{i,j}^{(n)}. \]
 Hence we conclude that $(\varphi_{n,0}|_{A_n})_* = \id _{G_n}$. 
 
{\noindent{\bf Proof\ of\ (1)}.\ }
Let $f\in C(\tX_n)\otimes M_{\td}$. When $\varphi_{n,0} (f) ^{-1} (\{ 0 \}) \cap (\tI_{0,j'}^{(n)}\cup\tI_{1,j'}^{(n)}) $ is $\varepsilon $-dense in $\tI_{0,j'}^{(n)}\cup\tI_{1,j'}^{(n)} $, by the definition of $D(t)$, we have that $f^{-1}(\{ 0 \})\cap (\tI_{0,j'}^{(n)}\cup\tI_{1,j'}^{(n)}) $ is $1/3\varepsilon$-dense in $\tI_{0,j'}^{(n)}\cup\tI_{1,j'}^{(n)} $. When $\varphi_{n,0}(f)^{-1} (\{ 0\}) \cap I_{i',j'}^{(n)}$ is $\varepsilon$-dense in $I_{i',j'}^{(n)}$ for any $(i',j')\in S^{(n)}$, by the definition of $\xi_{\delta}$, we have that $f^{-1} (\{ 0\}) \cap I_{i,j}^{(n)} $ is $1/2\varepsilon$-dense in $I_{i,j}^{(n)}$ for any $(i,j) \in S^{(n)}$. If $\varphi_{n,0} (f)^{-1} (\{ 0\}) \neq \phi$, then $f(c_n) =0$ by $\xi_c$. Inductively, define relatively prime natural numbers $\tp_{0,m}$, $\tp_{1,m}$, and the unital embeddings $\varphi_{n,m} : A_{n,m-1} \hookrightarrow A_{n,m}=A(G_n, \tp_{0,m}, \tp_{1,m}) $ by the same way as $\varphi_{n,0}$, then $\widetilde{\varphi}_{n,0} = \varphi_{n,m} \circ \varphi_{n,m-1}\circ \cdots \circ \varphi_{n,0}: A_n \hookrightarrow A_{n,m}$ satisfies (1) (and (0)). Thus we may claim that $\varphi_{n,0}$ satisfies (1).

{\noindent{\bf Proof\ of\ (2)}.\ }
Let $x_0$ and $y_0 \in \tX_n$ be such that $d(x_0, y_0) \leq \varepsilon$ and let $\heta_{j,i}: \tX_n \rightarrow \tX_n$ for $i=-1,0,1$ be the continuous  maps such that 
\[\eta_j (f) = f \circ \heta_{j,0} \oplus f\circ \heta_{j,1} \oplus f\circ \heta_{j,-1}.\]  By the definition of $\eta_j$, it follows that 
\[ d (\heta_{j,i} (x) , \heta_{j,i} (y)) \leq 1/3 d(x,y),\quad x, y \in \tX_n.\]
Then we have that  $|\Tr_{3r_0^{(n)}\td} (\eta (f) (x_0) - \eta(f) (y_0))|$
\begin{eqnarray} 
&\leq& 1/(3r_0^{(n)}) \sum_{j=1}^{r_0^{(n)}} \sum_{i=-1,0,1} | \Tr_{\td} (f\circ \heta_{j,i} (x_0) - f\circ \heta_{j,i} (y_0)) | \nonumber \\
&\leq& \max \{ |\Tr ( f(x) - f(y)) |; d(x,y) \leq 1/3\varepsilon \}.\nonumber
\end{eqnarray} 
Let $\hxi_{\delta}:\tX_n \rightarrow \tX_n$ be the continuous maps such that $\xi_{\delta} (f) = f \circ \hxi_{\delta}$. By the definition of $\xi_{\delta}$, it follows that 
\[ d(\hxi_{\delta} (x), \hxi_{\delta}(y)) \leq 1/2 d (x,y) ,\quad x,y \in \tX_n.\]
Then we have that $|\Tr_{\kappa\td} (\xi (f) (x_0) - \xi (f) (y_0))|$
\begin{eqnarray} 
 &\leq& 1/\kappa \sum_{\delta=1}^{\kappa} |\Tr_{\td} ( f\circ \hxi_{\delta} (x_0) - f\circ \hxi_{\delta} (y_0)) | \nonumber \\
& \leq&  \max  \{ | \Tr (f(x) - f(y) ) | ; \ d(x,y) \leq 1/2 \varepsilon \}. \nonumber 
\end{eqnarray}
Combining the estimation for $\eta$ and $\xi$, we have that 
\begin{eqnarray}
& &| \Tr_{\td'} ( \varphi_{n,0}(f)(x_0) - \varphi_{n,0} (f) (y_0) )|  \nonumber \\
& &= | \Tr  (( \eta \oplus \xi \oplus \xi_c) (f) (x_0) - (\eta \oplus \xi \oplus \xi_c) (f (y_0))) | \nonumber \\
& &\leq 3r_0^{(n)}/m | \Tr_{3r_0^{(n)}\td} (  \eta (f) (x_0)  -  \eta (f)(y_0)) | \nonumber \\ 
& &\quad  + \kappa /m  | \Tr_{\kappa \td} (\xi (f) (x_0) - \xi(f)(y_0)) | \nonumber \\
& &\leq \max \{ | \Tr_{\td} ( f(x) -f(y))| ; \ d(x,y) \leq 1/2 \varepsilon \}. \nonumber 
\end{eqnarray}
Thus we conclude (2).

\

{\noindent{\bf Proof\ of\ (3)}.\ }
Let $h_k^{(n)} \in C(\tX_n)_+^1,$ $n,$ $k\in\N$ be such that $h_k^{(n)} \searrow \chi(E_n),$ $k \rightarrow \infty$. 
Then we have that $\lim_{k\rightarrow \infty }\Tr_{\td'} (\varphi_{n,0} (h_k^{(n)} \otimes 1_{\td})(x))$
\[ = \lim_{k} \Tr ((\eta\oplus \xi) (h_k^{(n)} \otimes 1_{\td}) (x) \oplus 0_{(m-(3r_0^{(n)} + \kappa))\td}),\quad x \in \tX_n.\]
By the definition of $\eta$ and $\xi$ we have that 
\begin{eqnarray}
\lim_{k} \eta(h_k^{(n)} \otimes 1_{\td} ) &\leq& \chi(E_n) \otimes 1_{3r_0^{(n)}\td}, \nonumber \\
\lim_{k} \xi(h_k^{(n)} \otimes 1_{\td} ) &\leq& \chi(E_n) \otimes 1_{\kappa \td}. \nonumber 
\end{eqnarray}
Then, by $m \geq 2(3r_0^{(n)} + \kappa)$, we have that
\[ \Tr_{\td'} ( \varphi_{n,0}(\chi(E_n) \otimes 1_{\td})(x)) \leq 1/2 \Tr (\chi(E_n) \otimes 1_{\td'}(x)), \quad x \in \tX_n .\]
  
\begin{theorem}\label{maintheorem}
Let $G$ be a countable abelian group. Then there exists a unital simple $C^*$-algebra $A$ which has a unique tracial state, is expressed as the inductive limit $C^*$-algebra of an inductive sequence of dimension drop algebras, and satisfies $(K_0(A), [1_A])\cong (\Z, 1)$ and $K_1 (A) \cong G.$        
\end{theorem}
     
\bprf
 Let $(G_n)_{n\in \N}$ be an increasing sequence of finitely generated abelian groups and $\gamma_n$ connecting maps such that $\displaystyle G= \lim_{\longrightarrow} (G_n, \gamma_n).$ 
We may assume that $\gamma_n$ is injective.
By Proposition \ref{Prop2} and \ref{Prop3} we inductively obtain dimension drop algebras $A_{n,0}$ and $A_{n}$ for $G_n$ and unital embeddings $\psi_n: A_{n} \hookrightarrow A_{n+1,0}$ and $ \varphi_{n+1,0}: A_{n+1,0} \hookrightarrow A_{n+1} $ such that $ (\psi_n)_* \cong \gamma_n$, $(\varphi_{n+1,0})_* \cong \id_{G_{n+1}}$,  $\psi_n$ satisfies the conditions (0,1) and (0,2) in Proposition \ref{Prop2}, and $\varphi_{n+1,0} $ satisfies the conditions (1) and (2) in Proposition \ref{Prop3}. Set $\varphi_n = \varphi_{n+1,0} \circ \psi_n$.   
Because of the conditions (0,1), (0,2), (1), and (2), we see that $\varphi_n$ satisfies a similar conditions to (1) and (2) and satisfies $( \varphi_n)_* = \gamma_n$.
Let $A$ be the inductive limit $C^*$-algebra $\displaystyle \lim_{\longrightarrow} (A_n, \varphi_n).$ It is not so hard to show that $A$ is unital, projectionless, and satisfies $(K_0(A), [1_A]) \cong (\Z, 1)$ and $K_1(A) \cong G .$

By the condition (0,1) and (1) we shall show that $A$ is simple. 
Let $\widetilde{\varphi}_n$ be the canonical embedding from $A_n$ into $A$, and let $I$ be an ideal of $A$. Then it follows that 
\[ I =  \overline{(\bigcup \widetilde{\varphi}_n (A_n) \cap I)}.\]
Thus we obtain an ideal $I_n$ of $A_n$ such that $\widetilde{\varphi}_n (A_n) \cap I = \widetilde{\varphi}_n (I_n)$ and the compact subset $Y_n \subset X_n$ such that 
\[ I_n= \{ f \in A_n ;\ f|_{Y_n} =0\} .\]
Assume that $A\neq I$, then there exists a subsequence $m_n \in \N$ such that $Y_{m_n} \neq \phi,$ $n \in \N .$ For any $n_0\in \N$ and $f \in I_{n_0}$ we have that $\varphi_{m_n-1}\circ \varphi_{m_n-2} \circ \cdots \circ \varphi_{n_0} (f) \in I_{m_n},$ for large $n\in \N.$  By (0,1) and (1) of $\varphi_n$ it follows that $f^{-1} (\{ 0 \}) $ is $2^{-m}$-dense, where $m \rightarrow \infty$ $(m_n \rightarrow \infty)$, and then $f=0$, thus $I_{n_0} =\{ 0\}$ for any $n_0 \in \N$, which implies that $I= \{0\}$. Then we conclude that $A$ is simple.

By the conditions (0,2) and (2), for $f \in A_{n}$ and $\varepsilon >0 $ we obtain $\delta >0 $ and $M \in \N$ such that if $m \geq M$ then  
\begin{eqnarray}
& & \max \{ |\Tr(\varphi_{m,n}(f) (x)- \varphi_{m,n} (f)(y))| ;\ x,y  \in X_{m+1}\} \nonumber \\ 
&\leq & \max \{ |\Tr(\varphi_{m,n}(f) (x)- \varphi_{m,n} (f)(y))| ;\ x,y  \in \tX_{m+1},\ d(x,y)\leq 4 \pi\} \nonumber \\ 
&\leq & \max \{ |\Tr(\varphi_{m-1,n}(f) (x)- \varphi_{m-1,n} (f)(y))| ;\ x,y  \in \tX_{m},\ d(x,y)\leq 2 \pi\} \nonumber \\ 
&\vdots& \nonumber \\
&\leq &  \max \{ |\Tr_{} (f(x)-f(y))| ;\ x,y \in \tX_n,\ d(x,y) \leq  \delta \} \leq 1/2 \varepsilon,\nonumber
\end{eqnarray}
where $\varphi_{m,n}= \varphi_m\circ \varphi_{m-1} \circ \cdots \circ \varphi_n$.
Let $\tau_i \in T(A_{m+1})$, $i=0,1,$ and let $\mu_{\tau_i}$ be the probability measures on $X_{m+1}$ such that 
\[\tau_i \circ \varphi_{m,n} (f') = \int _{X_{m+1} } \Tr (\varphi_{m,n} (f') (x) ) d \mu_{\tau_i}(x),\quad  {\rm for\ any\ } f' \in A_n . \] By the above condition,  it follows that 
 $| \tau_0 \circ \varphi_{m,n} (f) - \tau_1 \circ \varphi_{m,n} (f) |$
\begin{eqnarray} 
&\leq& \sum_{i=0,1} |\tau_i \circ \varphi_{m,n} (f) - \Tr (\varphi_{m,n} (f) (x_0)) | \quad {\rm for\ some\ }x_0 \in X_{m+1} \nonumber \\
&\leq& \varepsilon, \nonumber 
\end{eqnarray}
for any $m \geq M.$ 
Then for a sequence $\tau_m \in T(A_m)$, $m\in \N$ it follows that  $(\tau_{m+1}\circ \varphi_{m,n} (f) )_{m \in \N,\ m \geq n} $ is a Cauchy sequence and converges to a point which is independent of the choice of $\tau_m \in T(A_m)$. We define a tracial state $\tau$ on $A$ by 
\[\tau (\widetilde{\varphi}_n (f) ) = \lim_{m \rightarrow \infty} \tau_{m+1} \circ \varphi_{m,n} (f) , \quad f\in A_n,\]
(remark that $\tau (\widetilde{\varphi}_{l+1}\circ \varphi_{l,n}(f)) = \tau(\widetilde{\varphi}_{n} (f))$). By the above claim for $\tau_0$ and $\tau_1$, it follows that $\tau$ is unique.  This completes  the proof. \eprf

\section{ $\Zj$-stability of projectionless $C^*$-algebras}\label{Sec4}
In this section, we prove that the projectionless $C^*$-algebra, that was constructed in Section \ref{Sec3}, for a countable abelian group without free part absorbs the Jiang-Su algebra tensorially without passing through the classification theory of these algebras (Theorem \ref{theorem2}). In order to prove $\Zj$-stability of the projectionless $C^*$-algebras for general countable abelian groups,  we shall modify the constructions in Section \ref{Sec3} ( Lemma \ref{lemma4.3} and Theorem \ref{maintheorem2}).

The following lemma is a continuation of the proof of Proposition \ref{Prop3}.  Remark that if $G_n=G_{n,1}$ then we have that $\tX_n=X_n$. We set $\natural (0,i)= i$ for $i=0,1$, ($\natural (i) = 1-i$, $d_{i,j}^{(n)} | \tp_{\natural (\natural(i,j))}$,) $\td=\tp_0 \tp_1$, $\td'= P_0P_1\td$, $A_n= A(G_n, \tp_0, \tp_1)$, and $A_{n,0}=A(G_n,P_0\tp_0, P_1\tp_1)$. Let $\varphi_{n,0}$ be the unital embedding from $A_n$ into $A_{n,0}$ that was defined in the proof of Proposition \ref{Prop3} and let $A_n$ and $A_{n,0}$ be defined by unital embeddings $\{ \epsilon_{i,j}^{(n)} \}$ and $\{ \epsilon _{i,j}^{(n,0)} \}$.  
\begin{lemma}\label{lemma4.1}
Suppose that $G_n=G_{n,1}$. Then there exists a unital embedding $\psi$ of $I(P_0, P_1)$ into $A_{n,0}$ satisfying that:
 if $f\in A_n$ and $\varepsilon >0$ with $\| f(x) - f(y) \| < \varepsilon$ for any $x$ and $y\in X_n$  then it follows that \[ \| [\psi(g), \varphi_{n,0} (f) ] \| < 2 \varepsilon,\quad {\ \rm for\ any \ } \ g \in I (P_0,P_1)^1.\] 
\end{lemma}
\bprf 
We define $C^*$-subalgebras $A_{i,j}$, $(i,j)\in S^{(n)}$ of $A_n$ by
\[A_{i,j} = \{ f\in A_n; f(x) = f(a_{i,j}^{(n)}) ,\ x \in  I_{i,j}^{(n)} \}.\quad \]
By the definition of $\varphi_{n,0}$ in the proof of Proposition \ref{Prop3}, it follows that 
\begin{eqnarray}
& &M_{\tp_i} \cong \varphi_{n,0} (A_{0,i}) (a_{0,i}^{(n)}) \subset \epsilon_{0,i}^{(n,0)} (M_{P_i\tp_i}),\ i=0,1,  \nonumber \\
& & M_{\td/d_{i,j}^{(n)}}\cong \varphi_{n,0} (A_{i,j}) (a_{i,j}^{(n)}) \subset \epsilon_{i,j}^{(n,0)}(M_{\td'/d_{i,j}^{(n)}}),\ (i,j)\in S^{(n)}\ {\rm with}\ i\geq 1  \nonumber
\end{eqnarray}
(where we denote by $A(x)$ the $C^*$-subalgebra $\{f(x); f \in A \}$ of $M_{\td'}$ for a $C^*$-subalgebra $A$ of $C(X_n)\otimes M_{\td'}$ and $x\in X_n$,) then there exist unital embeddings $\delta_{0,i} : M_{P_i} \hookrightarrow  \varphi_{n,0} (A_{0,i})(a_{0,i}^{(n)})' \cap \epsilon_{0,i}^{(n,0)}( M_{P_i \tp_i})$ for $i=0$,$1$, and $\delta_{i,j}:M_{P_{\natural (i,j)}} \hookrightarrow \varphi_{n,0} (A_{i,j}) (a_{i,j}^{(n)})' \cap \epsilon_{i,j}^{(n,0)} (M_{\td'/ d_{i,j}^{(n)}})$, for $(i,j) \in S^{(n)}$ with $i\geq 1$.

We denote by $\psi': C(X_n)\otimes M_{P_0 P_1} \hookrightarrow C(X_n) \otimes M_{\td'}$ the natural unital embedding  defined by $\psi'(f) (x) = 1_{\td} \otimes f(x)$. By the definition of $\varphi_{n,0}$,  we have that $\psi'(C(X_n) \otimes M_{P_0P_1}) \subset \varphi_{n,0} (1_{C(X_n)}\otimes M_{\td})'.$ Set $I=I(P_0, P_1)$ and identify $I$ with the $C^*$-subalgebra $\{ f \in C(X_n)\otimes M_{P_0P_1}; f(a_{0,0}^{(n)}) \in M_{P_0}\otimes 1_{P_1}, f(a_{0,1}^{(n)})\in 1_{P_0}\otimes M_{P_1}, f(\iota_{i,j}^{(n)} (t)) = f(\iota_{0,\natural(i,j)}^{(n)}(t))\ {\rm for}\ (i, j) \in S^{(n)},\ i\geq 1,\ {\rm and} \ t\in [0,2\pi] \}.$ Since $\psi'(I) (a_{0,i}^{(n)})\cong M_{P_i}$ and $\varphi_{n,0} (1_{C(\tX_n)} \otimes M_{\td})(a_{0,i}^{(n)}) \supset \varphi_{n,0}(A_{0,i})(a_{0,i}^{(n)})$, we obtain a unitary $w_{0,i} \in U(C(X_n)\otimes M_{\td'})$ such that $\Ad w_{0,i} \circ \psi' (I) (a_{0,i}^{(n)}) = \delta_{0,i} (M_{P_i})$, $w_{0,i} (x) \in \varphi_{n,0} (A_{0,i})(a_{0,i}^{(n)})'$ for any $x\in I_{0,i}^{(n)}$, and $w_{0,i}(c_n) = 1$.
Similarly, since $\psi'(I) (a_{i,j}^{(n)})\cong M_{P_{\natural(i,j)}}$ and $\varphi_{n,0} (1_{C(\tX_n)} \otimes M_{\td})(a_{i,j}^{(n)}) \supset \varphi_{n,0}(A_{i,j})(a_{i,j}^{(n)})$ for $(i,j)\in S^{(n)}$ with $i \geq 1$, we obtain a unitary $w_{i,j} \in U(C(X_n)\otimes M_{\td'})$ such that $\Ad w_{i,j} \circ \psi' (I) (a_{i,j}^{(n)}) = \delta_{i,j} (M_{P_{\natural(i,j)}})$, $w_{i,j} (x) \in \varphi_{n,0} (A_{i,j})(a_{i,j}^{(n)})'$ for any $x\in I_{i,j}^{(n)}$, and $w_{i,j}(c_n) = 1$.
We define $w \in U (C(\tX_n)\otimes M_{\td'}) $ by 
\[w(x)=w_{i,j}(x) ,\quad  x \in I_{i,j}^{(n)},\ (i,j)\in S^{(n)},\]
and $\psi = \Ad w \circ \psi '$, then it follows that $\psi (I) \subset A_{n,0}$.

Let $f\in A_n$ satisfy the condition in the lemma and let $f_{i,j} \in A_{i,j}$ for $(i,j)\in S^{(n)} $ be such that $ f_{i,j}(a_{i,j}^{(n)}) =f (a_{i,j}^{(n)}) $. By the definition of $\varphi_{n,0}$, we have that 
\[ \| \varphi_{n,0} (f_{i,j})(a_{i,j}^{(n)}) - \varphi_{n,0}(f) (x) \| = \|\bigoplus_{\delta=1}^{m} f_{i,j} (a_{i,j}^{(n)}) - f(x_{\delta}) \| < \varepsilon, \]
where $x_{\delta}\in X_n $, $m=P_0P_1$. 
Hence it follows that for $g \in I(P_0,P_1)^1$ and $x \in I_{i,j}^{(n)}$, $(i,j)\in S^{(n)}$, $\psi(g) \varphi_{n,0} (f) (x)$ $\approx_{\varepsilon} $ $w_{i,j} \psi' (g) w_{i,j}^* (x)$ $\varphi_{n,0} (f_{i,j}) (a_{i,j}^{(n)})
$ $= \varphi_{n,0} (f_{i,j}) (a_{i,j}^{(n)})$ $ w_{i,j} \psi' (g)$ $w_{i,j}^* (x) 
$ $\approx_{\varepsilon} \varphi_{n,0} (f) \psi (g) (x)$.
Whence $\| [\psi (g), \varphi_{n,0} (f)] \| < 2 \varepsilon$.
 \eprf

\begin{lemma}\label{lemma4.2}
Let $p$ and $q$ be relatively prime numbers. If relatively prime numbers $P_0$ and $P_1$ satisfy $\min \{P_0, P_1\} \geq 2pq$, then there exists a unital embedding $\psi$ of $I(p,q)$ into $I(P_0,P_1)$. 
\end{lemma}
\bprf Since $(p,q) = 1$ and $\min \{ P_0, P_1 \} > pq$, we obtain $a$, $b$, $c$, and $d\in \Z_+$ such that $P_0 = ap + bq$, $P_1=cp+dq$, $b < p$, and $c < q$. Set $D= ad -bc$. Since $\min \{ P_0, P_1 \} \geq 2 pq$ it follows that $D=$ $\frac{1}{pq} (P_0P_1 - P_1bq - P_0 cp)$ $ >\frac{1}{pq}((P_0-pq)(P_1-pq)-(pq)^2)$ $\geq0$. Then we can separate $P_0P_1$ as the following summation
\[ P_0P_1 = pcP_0+ qbP_1  + pqD .\]
We define a unital embedding $\psi: C( I_{0,0}^{(n)} \cup I_{0,1}^{(n)}) \otimes M_{pq} \hookrightarrow C(I_{0,0}^{(n)} \cup I_{0,1}^{(n)}) \otimes M_{P_0P_1}$ by 
\[\psi(f) (x) = V_0(f)\otimes 1_{cP_0}\oplus V_1(f)\otimes 1_{bP_1} \oplus f(x) \otimes 1_D,\] 
where $V_0$ and $V_1$ are the irreducible representations of $I(p,q)$ at $a_{0,0}^{(n)}$ and $a_{0,1}^{(n)}$. Since $cP_0 + qD =aP_1$ and $bP_1 + pD = d P_0$, we obtain unital embeddings $\epsilon_{0,i}:M_{P_i} \hookrightarrow M_{P_0P_1}$ such that $\psi(f)(a_{0,i}^{(n)}) \in \epsilon_{0,i}(M_{P_i})$, for any $f \in I(p,q)$ and $i=0,1$. Hence we have that $\psi(I(p,q)) \subset I(P_0, P_1)$, where $I(P_0,P_1)$ is defined by $\epsilon_{0,i}$. \eprf   

\begin{theorem}\label{theorem2}
The unital simple projectionless $C^*$-algebra, that is constructed in Theorem \ref{maintheorem}, for a countable abelian group without free part  absorbs the Jiang-Su algebra $\Zj$.
\end{theorem}
\bprf 
By the definition of $\varphi_{n,0}$ in the proof of \ref{Prop3},  for any $\varepsilon > 0$, $n_0\in \N$, and any finite subset $F \subset A_{n_0}$ there exists a natural number $m \in \N$ such that $\max $ $\{ \|\varphi_{m,n} (f) (x)$ $ - \varphi_{m,n}(f)(y) \|$ $ ;f\in F,\ x,y\in X_{m+1} \}$ $< \varepsilon$, where $\varphi_{m,n}= \varphi_m\circ \varphi_{m-1} \circ \cdots \circ \varphi_n$ and $\varphi_n =  \varphi_{n+1,0} \circ \psi_n: A_n \hookrightarrow A_{n+1}$ that was defined in the proof of Theorem \ref{maintheorem}. Applying Lemma \ref{lemma4.1} and Lemma \ref{lemma4.2}, we have that 
\[ I (p, q) \subset_{\rm unital} A' \cap A^{\infty}, \]
for any relatively prime  numbers $p$ and $q$.
By Proposition 2.2 in \cite{TW}, this implies that $A\otimes \Zj \cong A .$
\eprf

\begin{lemma}\label{lemma4.3}
Let $G_n$ be a finitely generated abelian group, $p$ and $q$ relatively prime numbers, $\tp_0$ and $\tp_1$ relatively prime numbers such that $d_{i,j}^{(n)} | \tp_0 \tp_1=\td$ for any $(i,j) \in S^{(n)}$ with $i \geq 1$, and $A_n$ the dimension drop algebra for $G_n$, $\tp_0$, and $\tp_1$. Then there exist relatively prime numbers $\tp_0'$ and $\tp_1'$ such that $\tp_0| \tp_0'$ and $\tp_1 | \tp_1'$, automorphisms $\{ \epsilon_j^{(n,0)} \in Aut (M_{\td'});$ $j=1,2,...,r_0^{(n)}\}$, where $\td'= \tp_0'\tp_1'$, unital embeddings $\{\epsilon_{i,j}^{(n,0)}:M_{\td'/d_{i,j}^{(n)}} \hookrightarrow M_{\td'};$ $(i,j)\in S^{(n)}\}$, where $d_{0,i}^{(n)} = \tp_{\natural(i)}'$ for $i=0,1$, and a unital embedding $\varphi_{n,0} :C(\tX_n)\otimes M_{\td} \hookrightarrow C(\tX_n) \otimes M_{\td'}$ satisfying the conditions (0) $\sim$ (3) in Proposition \ref{Prop3} and 
\begin{description}\label{prop:embedding}
\item[(4)]\label{(4)} 
there exist a unital embedding $\psi$ of $I(p,q)$ into $ A(G_n,\tp_0',\tp_1') $ satisfying that:  if $f\in A_n$ and $\varepsilon >0$ satisfy $\| f(x) - f(y) \| < \varepsilon$ for any $x$ and $y\in X_n$  then it follows that \[ \| [\psi(g), \varphi_{n,0} (f) ] \| < 2 \varepsilon,\quad {\ \rm for\ any \ } \ g \in I (p,q)^1.\] 
 \end{description} 

\end{lemma}
\bprf 
Let $\alpha$, $\beta \in \N$ be such that $\alpha p - \beta q =1$, $\alpha < q$, and $\beta <p$ and let $P_i$, $i=0,1$ be relatively prime numbers as in the proof of Proposition \ref{Prop3}, satisfying moreover $P_0P_1>$
$ (\alpha p + \beta q ) + \kappa $, $(P_i,pq)=1$,  and $\min \{P_0, P_1 \}\geq 2 pq $. Set $\tp_i' = P_i \tp_i$, $m=P_0P_1$, and $\td'=\tp_0'\tp_1'$.

Define a $*$-homomorphism $\eta: C(\tX_n)\otimes M_{\td} \rightarrow C(\tX_n)\otimes M_{(\alpha p+\beta q)\td}$  by 
\begin{eqnarray}
\eta(f )(x)=\left\{ \begin{array}{ll}
f(\tiota_{i,j}^{(n)} (t))\otimes 1_{\alpha p} \oplus f (\tiota_{\natural (i),j}^{(n)}(t))\otimes 1_{\beta q} ,\quad & x= \tiota_{i,j}^{(n)}(t),\ t \in [0,\pi], \nonumber \\
f(c_n)\otimes 1_{\alpha p + \beta q}, \quad & {\rm otherwise }, \nonumber 
\end{array} \right.\nonumber
\end{eqnarray} 
 define $\xi : C(\tX_n) \otimes M_{\td} \rightarrow C(\tX_n) \otimes M_{\kappa\td}$ and $\xi_c: C((\tX_n) \otimes M_{\td} \rightarrow C(\tX_n) \otimes M_{(m - (\alpha p +\beta q +\kappa))\td}$ by the same way as in the proof of Proposition \ref{Prop3}, and define $\varphi_{n,0}: C(\tX_n) \otimes M_{\td} \hookrightarrow C(\tX_n)\otimes M_{\td'}$ by $\varphi_{n,0}=$ $\eta \oplus$ $\xi \oplus$ $ \xi_c$. Remark that this $\varphi_{n,0}$ satisfies the condition (3).  Let $v' \in U(C(\tX_n)\otimes M_{\td})$ be such that 
$\Ad v' (f) (b_{1,j}^{(n)}) = f(b_{0,j}^{(n)})$ for any $f \in A_n$ and $j=1,2,...,r_0^{(n)}$, and $v' (x) =1_{\td}$ for any $x \in \bigcup_{(i,j)\in S^{(n)}} I_{i,j}^{(n)}$, and set 
\[v=v'\otimes 1_{\alpha p} \oplus v'^* \otimes 1_{\beta q} \oplus 1_{(m-(\alpha p + \beta q)) \td} \in U(C(\tX_n)\otimes M_{\td'}). \]
Thus it follows that $\Ad v \circ \varphi_{n,0} (f) (b_{1,j}^{(n)}) = \varphi_{n,0} (f) (b_{0,j}^{(n)})$ for any $f\in A_n$ and $j=1,2,...,r_0^{(n)}$. As in the proof of Proposition \ref{Prop3} for $\xi$, we obtain unital embeddings $\epsilon_{i,j}^{(n,0)}$, $(i,j)\in S^{(n)}$ such that  $\varphi_{n,0} (f) (a_{i,j}^{(n)}) \in \epsilon_{i,j}^{(n,0)} (M_{\td'/ d_{i,j}^{(n)}})$ for any $f \in A_n$ and $(i,j)\in S^{(n)}$. Let $A_{n,0}$ be the dimension drop algebra for $G_n$, $\tp_0'$, and $\tp_1'$ depending on the automorphisms $\{\Ad v(b_{1,j}^{(n)})\}$ and the unital embeddings $\{\epsilon_{i,j}^{(n,0)}\}$, thus it follows that 
$\varphi_{n,0} (A_n) \subset A_{n,0}$. Since $\alpha p - \beta q =1$, we have that $(\varphi_{n,0}|_{A_n})_* = \id_{G_n}$. 

As in the proof of Lemma \ref{lemma4.2} and  $(P_i,pq)=1$, we obtain $a$, $b$, $c$, $d$, and $D\in \N$ such that $D=ad -bc$ and $P_0P_1 =pcP_0 + q b P_1 + pqD $. Set $I=I(p,q)$, and define a unital embedding $\psi' : I \hookrightarrow C(\tX_n)\otimes M_{\td'}$ by $\psi' (g) (x) $
\[ = 1_{\td} \otimes (V_0 (g) \otimes 1_{\alpha} \oplus  V_1 (g) \otimes 1_{\beta} \oplus g (x)\otimes 1_D \oplus V_0(g) \otimes 1_{cP_0-\alpha}  \oplus V_1 (g) \otimes 1_{bP_1 - \beta}), \] 
where $V_0$ and $V_1$ are the representations of $I$ at $a_{0,0}^{(n)}$ and $a_{0,1}^{(n)}$, and we identify $I$ with the $C^*$-subalgebra $\{f\in C(\tX_n)\otimes M_{pq}; f(a_{0,0}^{(n)}) \in M_p \otimes 1_q , f(a_{0,1}^{(n)})\in 1_p \otimes M_q , f(\iota_{i,j}^{(n)} (t)) = f(\iota_{0,\natural(i,j)}^{(n)} (t))$ for $(i,j)\in S^{(n)}$ with $i \geq 1$ and $t \in [0,2 \pi]$, and $f(x) =f(c_n) $ for $x \in \tI_{i,j}^{(n)}$, $i=0,1$ and $j=1,2,...,r_0^{(n)} \}$. 
Then, for any $(i,j)\in S^{(n)}$, it follows that $\psi'(I) (a_{i,j}^{(n)})$ is contained unitally in a $C^*$-subalgebra of $\varphi_{n,0}(1_{C(\tX_n)}\otimes M_{\td}) (a_{i,j}^{(n)})'$ $\cap M_{\td'}$ that is isomorphic to $M_{P_{\natural(i,j)}}$. As in the proof of Lemma \ref{lemma4.1}, we obtain a unitary $w\in U(C(\tX_n)\otimes M_{\td'}$ such that $w(x) =1_{\td'}$ for $x \in \tI_{i,j}^{(n)}$, $\Ad w \circ \psi' (g) (a_{i,j}^{(n)}) \in \epsilon_{i,j}^{(n,0)} (M_{\td/d_{i,j}^{(n)}})$ for $g \in I$ and $(i,j)\in S^{(n)}$,  and 
\[\| [ \Ad w \circ \psi' (g) (x), \varphi_{n,0} (f) (x) ] \| < 2 \varepsilon,\]
for $g \in I^1$, $x\in \tX_n$, $f\in A_n $, and $\varepsilon >0 $ satisfying that $\| f(x) - f(y) \| < \varepsilon $ for any $x, y\in \tX_n$.
Set $\psi = \Ad w \circ \psi'$, then by the definitions of $v$ and $\psi'$, it follows that 
\[ \Ad v \circ \psi (g) (x) = \psi (g) (x) , \quad x \in \tX_n,\]
which induces that $\psi (I) \subset A_{n,0}$.

In order to obtain the condition (1) and (2)  for this $\varphi_{n,0}$, applying Proposition \ref{Prop3},  we take the dimension drop algebra $A_{n,1}=$ $A(G_n,\tp_0'', \tp_1'')$ and the unital embedding $\varphi_{n,1}:A_{n,0} \hookrightarrow A_{n,1}$, then $\varphi_{n,1}\circ \varphi_{n,0}$ satisfies (1) and (2). Hence we can construct $\varphi_{n,0}$ satisfies (0) $\sim$ (4). 
\eprf 
 
Applying Proposition \ref{Prop2} and Lemma \ref{lemma4.3} in place of Proposition \ref{Prop3}, we conclude the following theorem.

\begin{theorem}\label{maintheorem2}
Let $G$ be a countable abelian group. Then there exists a unital simple projectionless $C^*$-algebra $A$ which has a unique tracial state, is expressed as the inductive limit $C^*$-algebra of an inductive sequence of dimension drop algebras, satisfies $(K_0(A), [1_A])\cong (\Z, 1)$, $K_1 (A) \cong G$, and absorbs the Jiang-Su algebra $\Zj$.        
\end{theorem}

\section{ Certain aperiodic automorphism of projectionless $C^*$-algebras }\label{Sec5}

In this section, we show the UHF-embeddability of the $C^*$-algebra in Theorem \ref{maintheorem} or \ref{maintheorem2} ( Proposition \ref{prop:UHF}) and prove  Theorem \ref{maintheorem3}.

Let $A$ be the unital simple projectionless $C^*$-algebra which is obtained in Theorem \ref{maintheorem} or \ref{maintheorem2} for a countable abelian group $G$, $\tau$ the unique tracial state of $A$, $(A_n, \varphi_n)$ the sequence of dimension drop algebras such that $\displaystyle A= \lim_{\longrightarrow} (A_n, \varphi_n)$, where we have defined $\varphi_n:A_n \hookrightarrow A_{n+1}$ as the decomposition of unital embeddings such that
\[ A_n =A(G_n, \tp_{0,n}, \tp_{1,n}) \longrightarrow \hspace{-1.50em}\raisebox{1.ex}{$^{\psi_n}$}\ \ A(G_{n+1}, \tp_{0,n+1}' , \tp_{1,n+1}')=A_{n+1,0} \longrightarrow \hspace{-1.75em}\raisebox{1.ex}{$^{\varphi_{n+1,0}}$}A_{n+1}.\] Set $\td_n= \tp_{0,n}\tp_{1,n}$ and $\td_{n,0}=\tp_{0,n}'\tp_{1,n}'$.
 
Let $\mu_n$ be the probability measure on $X_n$ defined by $\tau \circ \tphi_n ( f\otimes 1_{\td_n} ) $ $= \int_{X_n} f(x)$ $ d\mu_n (x)$ for any $f\in C(X_n)$, and $E_n = \{ a_{i,j}^{(n)} \} \cup \{\tau_j^{(n)} (z(\pi)) \} \subset X_n $.
We recall that $\psi_n$  and $\varphi_{n+1,0}$ satisfy (0,3) and (3) in Proposition \ref{Prop2} and \ref{Prop3} or Lemma \ref{lemma4.3}, then we have the following lemma. 
\begin{lemma}\label{lem:uniq.tr}
For any $n \in \N$ we have that $\mu_n (E_n) =0$. 
 \end{lemma}
\bprf
Let $h_{k}^{(n)} \in C (X_n)_+^1$ be such that $h_{k}^{(n)} \searrow \chi (E_n)$ (pointwise). By the conditions (0,3) and (3) in Proposition \ref{Prop2} and \ref{Prop3} or Lemma \ref{lemma4.3} we have that 
\begin{eqnarray}
\mu_n (E_{n}) &=& \lim_{k \rightarrow \infty} \tau \circ \tphi_n (h_k^{(n)} \otimes 1_{\td} ) = \lim_k \tau \circ \tphi_{n+1} \circ \varphi_{n+1,0} \circ \psi_n  (h_k^{(n)}\otimes 1_{\td}), \nonumber \\
&\leq & \lim \tau \circ \tphi_{n+1} \circ \varphi_{n+1,0} (h_k^{(n+1)} \otimes 1_{\td_{n+1,0}}) \nonumber \\
& \leq & \lim \int_{X_{n+1}} \Tr_{\td_{n+1}} (\varphi_{n+1,0} (h_k^{(n+1)} \otimes 1_{\td_{n+1,0}})(x)) d \mu_{n+1} (x) \nonumber \\
& \leq & 1/2 \lim \int_{X_{n+1}} \Tr_{\td_{n+1}} (h_k^{(n+1)} \otimes 1_{\td_{n+1}}(x)) d \mu_{n+1} (x) \nonumber \\
&= &1/2 \mu_{n+1} (E_{n+1}). \nonumber  
\end{eqnarray}
Hence we have that $\mu_n(E_n) =0$.
\eprf

Let $(\pi_{\tau}, H_{\tau})$ be the GNS-representation of $A$ associated with $\tau$.
The following Proposition is a generalization of Proposition 2.3 in \cite{Sat}.\begin{proposition}\label{prop:UHF}
 There exists a UHF $C^*$-subalgebra $B$ of $\pi_{\tau} (A) ''$ and an increasing sequence $(B_n)_{n\in \N}$ of matrix $C^*$-subalgebras of $B$ such that 
$\pi_{\tau} (A) \subset B$, $1_{B_n} = 1_{B}$, $B= \overline{\bigcup B_n}^{\| \cdot \|}$, and 
\[ ( B_n' \cap \pi_{\tau}(A) )'' = B_n'\cap \pi_{\tau} (A) ''.\] 
\end{proposition}
\bprf Let $C_n$ denote the $C^*$-algebra $C (\tX_n) \otimes M_{\td_n}$, where $\td_n$ is such that $A_n \subset C(\tX_n)\otimes M_{\td_n}$ and let $f_k \in C(X_n )_+^1,$ $k\in \N$ be such that $f_k \nearrow \chi(X_{n}\setminus E_n)$, $k\rightarrow \infty$ (pointwise), and $f_k(E_n)=\{ 0\}$ and define a map $\Psi_n : C_n \rightarrow \pi_{\tau}(A)''$ by 
\[ \Psi_n (g) = \lim_{k \rightarrow \infty} \pi_{\tau}\circ \tphi_n((f_k \otimes 1_{\td}) g), \quad g \in C_n,\]
where the limit is taken in the strong operator topology.
From Lemma \ref{lem:uniq.tr} and a similar argument to the proof of Proposition 2.3 in \cite{Sat},
it follows that $\Psi_n$ is independent of the choice of $f_k$, is an injective $*$-homomorphism, $\Psi_n|_{A_n} = \pi_{\tau}\circ \tphi_n$ and satisfies $\Psi_n =\Psi_{n+1} \circ \varphi_n$ on $C_n$. Then we have that 
\[ \pi_{\tau}(A) \subset \overline{\bigcup_{n} \Psi_n(C_n)}^{\| \cdot \|} \subset \pi_{\tau}(A)''.\]
Set $B = \overline{\bigcup_{n} \Psi_n(C_n)}^{\| \cdot \|}$ and define a sequence of finite dimensional $C^*$-subalgebras $B_n$ of $B$ by 
\[B_n = \Psi_n (1_{C(\tX_n)}\otimes M_{\td_n}),\quad n\in \N .\]
By the definition of $\varphi_{n+1,0}$ and $\psi_n$, for $f\in C_n$ there exists a large natural number $m$ such that $\varphi_m \circ \varphi_{m-1} \circ \cdots \circ \varphi_n (f)$ is almost contained in $1_{C(\tX_{m+1})}\otimes M_{d_{m+1}}$, then we have that $\Psi_n (C_n) \subset \overline{\bigcup B_n},$ hence $B=\overline{\bigcup B_n}$ is a UHF-algebra.

We shall show $(B_n' \cap \pi_{\tau}(A) )'' \supset B_n' \cap \pi_{\tau}(A)'' $.
From Lemma \ref{lem:uniq.tr} we have that 
\[ \Psi_N ((1-f_k)\otimes 1_{\td_N}) = \pi_{\tau} \circ \tphi_N ((1-f_k)\otimes 1_{\td_N}) \searrow 0, \quad k\rightarrow \infty\ {\rm (strongly)} , \]
for any $N\in \N$. Let $g \in B_n' \cap \pi_{\tau} (A) '' $. We may assume that $g \in B_n' \cap \Psi_N(C_N) $, $N>n$. 
Set $B_{n,0} = \varphi_{N-1,n} (1\otimes M_{\td_n}) \subset C_N$.
We obtain $g_0 \in C_N \cap B_{n,0}'$ such that $g = \Psi_N (g_0).$ Since 
$g_0\cdot (f_k \otimes 1_{\td_N}) \in A_N \cap B_{n,0}'$ and $\Psi_N (g_0\cdot (f_k \otimes 1_{\td_N})) \rightarrow \Psi_N(g_0)$ (strongly), we have that $g\in (\pi_{\tau}(A) \cap B_n')''$.
 \eprf

\begin{remark}\label{remarkcor}
By the UHF-embeddability in the above proposition  and by the same argument as in  the proof of Theorem 2.6 in \cite{Sat}, we can prove \ref{maincorollary} for the projectionless $C^*$-algebras that are constructed in this paper.  

\end{remark}

The following lemma is an adaptation of Theorem 3.12.14 in \cite{Ped}.
\begin{lemma}\label{lem:autfundamental}
Suppose that a sequence $(F_n)_{n\in \N}$ of projections in $\pi_{\tau} (A) ''$  satisfies  that $[F_n,x] \rightarrow 0\ ({\rm strongly})$ for any $\ x \in \pi_{\tau} (A) '' $.   
Then there exists a central sequence $(f_n)_n $  of positive elements in $A^1$ (in the sense that 
\[ \| [f_n, a ] \|  \rightarrow 0 ,\ {\rm for \ any\ } a \in A,) \]
such that 
\[F_n - \pi_{\tau} (f_n)  \rightarrow 0\ \quad  ({\rm strongly })  .\]

\end{lemma}

We recall that \[ \WInn (A) = \{ \alpha \in \Aut (A) :\pi_{\tau}\circ \alpha= \Ad W \circ \pi_{\tau} ,\quad W\in U(\pi_{\tau}(A)'')\}.\]
Mimicking the proof of Lemma 4.4 and Theorem 4.5 in \cite{Kis}, we obtain the following proposition.
\begin{proposition}\label{prop:wouter}
Let $A$ be a unital simple projectionless $C^*$-algebra with a unique tracial state $\tau$ and let $\alpha \in \Aut (A)$. If there exists a central sequence $(f_n)_{n\in \N}\in (A_{\infty})_+^1$ such that 
\[ (\alpha (f_n))_n \cdot (f_n)_n =0,\quad   \tau (1- (f_n + \alpha (f_n))) \rightarrow 0,\ n\rightarrow \infty, \]
then $\alpha \notin \WInn (A)$ .
\end{proposition}
\bprf
Assume that there exists $V \in U(\pi_{\tau} (A) '')$ such that $\pi_{\tau} \circ \alpha = \Ad V \circ \pi_{\tau}.$ Let $\rho: A \times_{\alpha} \Z \rightarrow \pi_{\tau} (A)''$ be a representation of $A \times_{\alpha} \Z$ determined by $\rho (a) = \pi_{\tau} (a)$ for $a \in A$ and $\rho(u_{\alpha}) =V$, where $u_{\alpha}$ is the canonical unitary with $u_{\alpha}a = \alpha(a)u_{\alpha},$ $a \in A$. Define a tracial state $\phi$ of $A\times_{\alpha} \Z$ by 
$\phi(x) = (\rho(x) \Omega_{\tau} | \Omega_{\tau}),$ where $\Omega_{\tau}$ is the cyclic vector obtained by the GNS-representation associated with $\tau$. By $V \in U(\pi_{\tau} (A)'')$ there are $a_n \in A$ such that $\pi_{\tau} (a_n) \rightarrow V^*$ (strongly). Thus we have $\phi(a_n u_{\alpha}) \rightarrow 1$. However we shall obtain $\phi (au_{\alpha}) =0$ for any $a \in A$ by the condition of $\alpha$ and the following argument.  

Set $a \in A^1$ and let $(f_n) \in (A_{\infty})_+^1$ be as in the Proposition, then it follows that  
\begin{eqnarray} 
\phi (au_{\alpha}) &=& \lim_{n\rightarrow \infty} \phi ((f_n + \alpha(f_n))au_{\alpha})+ \phi((1-(f_n +\alpha (f_n)))a u_{\alpha}) \nonumber \\
&=& \lim \phi ((f_n + \alpha (f_n))au_{\alpha}). \nonumber 
\end{eqnarray}
Since $(f_n)_n$ is a central sequence we have that 
\[ \lim_{n\rightarrow \infty } \phi(f_na u_{\alpha}) = \lim \phi (f_n^{1/2}\alpha (f_n)^{1/2} a u_{\alpha}) =0 , \]
which implies that $\phi(au_{\alpha}) =0$. \eprf

The following proof of Theorem \ref{maintheorem3} is based on the proof of Theorem 4.5 in \cite{Kis}. 

\bprfm 
Suppose that $\alpha \in \Aut (A) $ has a central sequence $(f_n')_n \in (A_{\infty})_+^1$ which satisfies the condition in the theorem for $2k$ in place of $k$. Set $f_n = \sum _{j=0}^{k-1} \alpha^j (f_n').$ Thus $(f_n)_n \in (A_{\infty})_+^1$ satisfies 
$(f_n) \cdot \alpha^{k}((f_n)) =0,$ and $\tau(1_A - (f_n + \alpha ^k (f_n))) \rightarrow 0 $. From Proposition \ref{prop:wouter}, 
it follows that $ \alpha^k \notin \WInn (A)$. Then we have that  $[\alpha]$ is aperiodic in $\Aut (A) / \WInn (A) $.

Suppose that $[\alpha] \in \Aut(A) / \WInn (A) $ is aperiodic. Let $\oalpha$ be the weak extension of $\pi_{\tau} \circ \alpha \circ \pi_{\tau}^{-1} $ to an automorphism of $\pi_{\tau} (A) ''.$ 
Let   $k\in \N,$ $\varepsilon_n > 0$, $n\in \N$, be such that $\varepsilon_n \searrow 0,$ and $F_n$, $ n\in \N$, be finite subsets of $A^1$ such that $\overline{(\bigcup F_n)}^{\|\cdot \|} = A^1.$
By the classification theory for aperiodic automorphisms on the injective type $I\hspace{-.1em}I_1$ factor due to Connes \cite{Con}, there exists a central sequence $(E_m) \in (\pi_{\tau} (A) '')_{\infty}$  in the sense of the strong topology of projections such that 
\[E_m + \oalpha(E_m)+ \cdots + \oalpha^{k-1}(E_m) \rightarrow  1\ ({\rm strongly }).\]

By Lemma \ref{lem:autfundamental} we obtain a central sequence $(f_m')_m \in A_{\infty}$ in the sense of the norm such that $f_m'\in A_+^1$ and $E_m -\pi_{\tau} (f_m') \rightarrow  0 \ ({\rm strongly})$. Then we have that $\displaystyle \tau (1 - \sum_{j=0}^{k-1} \overline{ \alpha}^j \circ \pi_{\tau} (f_m')) \rightarrow 0, $ $m\rightarrow \infty$. 
Set 
\[ g_m = (f_m')^{1/2} ( \alpha (f_m') + \alpha^2 (f_m') + \cdots + \alpha^{k-1} (f_m')) (f_m') ^{1/2},\quad m\in \N, \]
then we have $\| \pi_{\tau} (g_m) - E_m ( \oalpha (E_m) + \oalpha^2 (E_m) + \cdots + \oalpha^{k-1} (E_m))E_m \|_2 $
\begin{eqnarray}
 &\leq & \| (\pi_{\tau} (f_m')^{1/2} - E_m)( \sum_{j=1}^{k-1} \oalpha^j (\pi_{\tau} (f_m') ))\pi_{\tau}(f_m')^{1/2} \|_2 \nonumber \\
 & + & \| E_m( \sum \oalpha^j (\pi_{\tau} (f_m') - E_m)) \pi_{\tau}(f_m')^{1/2}\|_2 \nonumber \\
 & + & \| E_m ( \sum \oalpha^j (E_m ))(\pi_{\tau} (f_m')^{1/2} - E_m) \|_2\nonumber \\
 & \leq & 2k \| \pi_{\tau} (f_m')^{1/2} - E_m \|_2 + k\| \pi_{\tau} (f_m') - E_m\|_2 \rightarrow 0 \nonumber 
\end{eqnarray} 
 and by $ \|E_m( \oalpha (E_m) + \oalpha^2 (E_m) + \cdots + \oalpha^{k-1} (E_m))E_m \|_2 \rightarrow 0,$ we have that $\| \pi_{\tau} (g_m) \|_2 \rightarrow 0.$

Define $f_{\varepsilon} \in C( [0, \infty))_+^1$ by 
\begin{eqnarray}
f_{\varepsilon} (t) =\left\{ \begin{array}{ll}
(2\varepsilon)^{-2} t ,\quad & 0 \leq t \leq 4 \varepsilon^2 , \\
1 \quad  &4 \varepsilon^2 < t,  
\end{array} \right.\nonumber
\end{eqnarray} 
and let
\[ f_{m,n} = (f_m')^{1/2} \cdot (1_A- f_{\varepsilon_n}(g_m)) \cdot (f_m')^{1/2}\in A_+^1 ,\quad m,n \in \N .\]
Remark that $0 \leq f_{m,n} \leq f_m'$ and $(f_{m,n})_m$ is a central sequence in the norm sense for any $n \in \N$ (i.e., $\|[f_{m,n}, a]\| \rightarrow 0,\ m \rightarrow \infty ,\ a \in A$).
Then for any $j=1,2,...,k-1$ and $m,$ $n \in \N$ we have that $\| f_{m,n} \alpha^j (f_{m,n})\|^2 $
\begin{eqnarray}
&\leq & \| f_{m,n} \alpha^j (f_{m,n})^{1/2} \|^2 = \| f_{m,n} \alpha^j (f_{m,n}) f_{m,n} \| \nonumber \\
& \leq & \|f_{m,n} (\alpha (f_m') + \alpha^2 (f_m') + \cdots + \alpha^{k-1} (f_m')) f_{m,n} \| \nonumber \\
& \leq & \| (f_m')^{1/2} (\alpha (f_m') + \cdots + \alpha ^{k-1} (f_m'))(f_m')^{1/2} (1_A- f_{\varepsilon_n}(g_m)) (f_m')^{1/2}\| \nonumber \\
& \leq & \| g_m  -g_m \cdot f_{\varepsilon_n} (g_m) \| = 
\| g_m \cdot (1- f_{\varepsilon_n}) (g_m) \| = \varepsilon_n^2. \nonumber
\end{eqnarray}

Since $\| \pi_{\tau}(g_m) \|_2 \rightarrow 0 $ and $(f_{m,n})_m \in A_{\infty}$, for $\varepsilon_n > 0$ and $F_n$ we obtain a sub-sequence $m_n \in \N$ such that \[ \|f_{m_n,n} - f_{m_n}' \|_2 \leq \| f_{\varepsilon_n} (g_{m_n}) \|_2 < \varepsilon_n , \]
\[ \| [ f_{m_n ,n} , x ] \| < \varepsilon_n , \quad x \in F_n .\]
Let $f_n = f_{m_n, n} \in (A)_+^1$.
Thus $(f_n)_n$ is a central sequence in the norm sense which satisfies 
\[ (f_n)_n (\alpha^j ( f_n ))_n = 0 , \quad j=1,2,...,k-1. \]
Since $\displaystyle \tau (1 - \sum_{j=0}^{k-1} \overline{ \alpha}^j \circ \pi_{\tau} (f_{m_n}')) \rightarrow 0, $ $(n\rightarrow \infty)$, we have $\displaystyle \tau (1 - \sum_{j=0}^{k-1} \alpha^j (f_n)) \rightarrow 0 .$ This completes the proof. 
\eprf

{\bf Acknowledgment.} 
The author is grateful to Professor Akitaka Kishimoto for his support and valuable comments.


\begin{thebibliography}{Pi1}
\bibitem[1]{Con}A. Connes. {\it Outer conjugacy class of automorphisms of factors}, Ann. Scient. Ec. Norm. Sup. (4) {\bf 8} (1975), 383-420.
\bibitem[2]{DW}M. Dadarlat and W. Winter {\it Trivialization of $C(X)$-algebras with strongly self-absorbing fibres,} Bull. Soc. Math. France 136 (2008), no. 4, 575-606.
\bibitem[3]{DPT}M. Dadarlat, N. C. Phillips, and A. S. Toms {\it A direct proof of $\Zj$-stability for AH algebras of bounded topological dimension, } arXiv:0806.2855 17 Jun 2008.

\bibitem[4]{Ell} G. A. Elliott {\it An invariant for simple $C^*$-algebras,}\ Canadian Mathematical Society. 1945-1995
3 (1996), 61-90.
\bibitem[5]{ET} G. A. Elliott and A. S. Toms {\it Regularity properties in the classification program for separable amenable $C^*$-algebras, } Bull. Amer. Math. Soc. 45 (2008), 229-245.

\bibitem[6]{JS} X. Jiang and H. Su. {\it On a simple unital projectionless $C\sp *$-algebra,}  Amer. J. Math.  121  (1999),  no. 2, 359--413.
\bibitem[7]{Kis}A. Kishimoto. {\it The Rohlin property for shifts on UHF algebras and automorphisms of Cuntz algebras,}  J. Funct. Anal.  140  (1996),  no. 1, 100--123. 
\bibitem[8]{Li} L. Li. {\it Classification of simple $C^*$-algebras: inductive limits of matrix algebras over trees,} Mem. Amer. Math. Soc. 127 (1997), no. 605, vii+123 pp.
\bibitem[9]{Lin}H. Lin. {\it Inductive Limits of Subhomogeneous $C^*$-algebras with Hausdorff Spectrum, } arXiv:0809.5273 30 Sep 2008.
\bibitem[10]{My} J. Mygind. { \it Classification of certain simple $C^*$-algebras with torsion in $K_1$,} Canad. J. Math. 53 (2001), no. 6, 1223-1308.
\bibitem[11]{Ped}G. K. Pedersen. ``$C\sp{*} $-algebras and their automorphism groups.'' London Mathematical Society Monographs, 14. Academic Press, Inc. [Harcourt Brace Jovanovich, Publishers], London-New York, 1979. 
\bibitem[12]{RW}M. R$\o$rdam and W. Winter. {\it The Jiang-Su algebra revisited,} arXiv:0801.2259 
\bibitem[13]{Sat}Y. Sato. {\it Certain aperiodic automorphisms of unital simple
projectionless $C^*$-algebras, } submitted. 
\bibitem[14]{Thomsen} K. Thomsen {\it Limits of certain subhomogeneous $C^*$-algebras,} Mem. Soc. Math. Fr. (N.S.) No. 71 (1997), vi+125 pp. (1998).
\bibitem[15]{TW}A. S. Toms and W. Winter. {\it Z-stable ASH algebras,}arXiv:math/0508218 
\bibitem[16]{W}W. Winter. {\it Localizing the Elliott conjecture at strongly self-absorbing $C^*$-algebras -with an appendix by Huaxin Lin-, } arXiv:0708.0283  11 Sep 2007.
\end{thebibliography}
\end{document}